 \documentclass[final,authoryear,3p,times,fleqn]{elsarticle}

\usepackage{amssymb}
\usepackage{amsmath}
\usepackage{amsthm}

\usepackage[ruled]{algorithm2e}

\usepackage{graphicx}
\usepackage{float} 
\usepackage{subcaption}
\renewcommand*{\thesubfigure}{(\alph{subfigure})}

\usepackage{pxfonts,times}
\usepackage{siunitx}

\usepackage{bm}
\usepackage{color}
\usepackage{comment}
\usepackage{natbib}

\usepackage[
	colorlinks=true, 
	citecolor=blue, 
	urlcolor=blue]{hyperref}

\newtheorem{theo}{Theorem}
\newtheorem{defi}{Definition}

\newtheorem{prop}{Proposition}

\newtheorem{coro}{Corollary}

\newtheorem*{prf}{Proof}

\makeatletter
\def\th@plain{\upshape}
\makeatother

\makeatletter
	
	\@addtoreset{equation}{section}
\makeatother

\journal{Elsevier}

\begin{document}

\begin{frontmatter}



\title{Dynamic system optimal traffic assignment with atomic users: \\ Convergence and stability}



\author[a]{Koki Satsukawa\corref{cor1}}
\ead{satsukawa@tohoku.ac.jp} 
\cortext[cor1]{Corresponding author. Tel.: +81-22-795-7492 (ext. 7492).}

\author[b]{Kentaro Wada}
\ead{wadaken@sk.tsukuba.ac.jp}
\author[c]{David Watling}
\ead{D.P.Watling@its.leeds.ac.uk}

\address[a]{New Industry Creation Hatchery Center, Tohoku University, Miyagi, Japan}
\address[b]{Faculty of Engineering, Information and Systems, University of Tsukuba, Ibaraki, Japan}
\address[c]{Institute for Transport Studies, University of Leeds, Leeds, United Kingdom}

\begin{abstract}
In this study, we analyse the convergence and stability of dynamic system optimal (DSO) traffic assignment with fixed departure times.
We first formulate the DSO traffic assignment problem as a strategic game wherein atomic users select routes that minimise their marginal social costs, called a `DSO game'.
By utilising the fact that the DSO game is a potential game, we prove that a globally optimal state is \textcolor{black}{a stochastically stable state} under the logit response dynamics, and the better/best response dynamics converges to \textcolor{black}{a locally optimal state}.
Furthermore, as an application of DSO assignment, we examine characteristics of the evolutionary implementation scheme of \textcolor{black}{marginal cost pricing}.
Through theoretical comparison with a fixed pricing scheme, we found the following properties of the evolutionary implementation scheme: (i) the total travel time decreases smoother to an efficient traffic state as congestion externalities are perfectly internalised; (ii) a traffic state would reach a more efficient state as the globally optimal state is stabilised.
Numerical experiments also suggest that these properties make the evolutionary scheme robust in the sense that they prevent a traffic state from going to worse traffic states with high total travel times.

\end{abstract}


\end{frontmatter}


\section{Introduction}\label{Sec:Introduction}
Dynamic system optimal (DSO) traffic assignment represents normative traffic flow patterns minimising \textcolor{black}{total costs} in transport networks. 
The optimal solutions of a DSO problem provide useful insights into the design of efficient transport management and control schemes, while the value of the objective function is the benchmark for evaluating these schemes.
Therefore, DSO assignment has attracted significant attention over the decades since the pioneering work of \cite{Merchant1978,Merchant1978a}.

There are two major research streams for analysing DSO assignment. 
The first analyses the properties of the mathematical optimisation problems of DSO assignment.
In this stream, several studies investigated the relationship between the solutions of total cost minimisation problems and route marginal social costs.
Regarding DSO assignment with fixed departure times, \cite{Carey1993} and \cite{Nie2011} showed that \textcolor{black}{route choice equilibrium} of marginal social costs is a first-order necessary condition for optimality in networks with many-to-one origin--destination (OD) pairs.
\cite{Ziliaskopoulos2000} and \cite{Carey2012} showed a convex approximation of DSO problems; they also demonstrated that marginal cost equilibrium is a sufficient condition for optimality in the problems\footnote{Although this approximation allows the existence of `vehicle-holding' situations, some studies developed methods for resolving this problem computationally~\citep[e.g.][]{Zhu2013,Shen2014}.}.
The second stream provides qualitative insights into the optimal control rules from simple networks, where exact optimal solutions can be analysed.
For example, \cite{Kuwahara2001}, \cite{Munoz2006} and \cite{Zhao2018} developed graphical solution methods in simple parallel-link networks based on the concept of marginal cost equilibrium, and derived properties of optimal ramp metering control.
\cite{Shen2009} and \cite{Zhang2010} investigated an optimal tolling scheme in a corridor network and a ramp control policy in a monocentric network, respectively.

Although the existing studies indicate some useful properties of the optimal states and costs, little is known about theoretical results on the convergence of dynamical processes (i.e. evolutionary dynamics and iterative algorithms) to optimal states and the stability of them, despite their importance in achieving the optimal states. 
For instance, several studies proposed heuristic solution algorithms~\citep[e.g.][]{Ghali1995,Shen2007,Qian2012,Zhang2020}; however, convergence is not guaranteed in most cases. 
The main reason for the difficulty in addressing the theoretical issues is associated with the non-convexity of DSO problems (or non-monotonicity of the mapping of variational inequality (VI) formulation of the marginal cost equilibrium problems), which stems from a complex dynamic loading model~\citep[][]{Carey1992}: as the first-order necessary condition is not sufficient for optimality, it is essentially difficult to establish the convergence to the globally optimal state of the standard dynamical processes (e.g. deterministic evolutionary dynamics and gradient based algorithms).
Furthermore, such dynamical processes may not converge to even a \textit{locally} optimal state of total cost minimisation problems (i.e. DSO solution in a broad sense)\footnote{One exception is \cite{Garcia2000}, who reported that the evolutionary process of the historical frequency of route usage induced by fictitious play converges to a local optimal state.
However, this does not mean that traffic flow patterns converge to a particular state.
Moreover, the type of local optimal states that are achieved through this process are not known; the resulting local optimal state may be worse than an initial state.}, since traffic states satisfying the first-order condition include local maxima and saddle points.

This study analyses convergence and stability properties of DSO assignment minimising the total travel time with fixed departure times in a general network.
\textcolor{black}{To this end, we employ a game theory approach that deals with atomic users and formulate dynamic traffic assignment problems in the framework of game theory.
This approach enables us to establish the relationship between the traffic assignment problems and certain classes of games where convergence and stability properties are well-known.}
Specifically, we first formulate the DSO traffic assignment problem as a strategic game wherein atomic users (individual vehicles) select routes that minimise their marginal social costs, called a `DSO game’.
By utilising the fact that the DSO game is a potential game~\citep[][]{Monderer1996}, we rigorously analyse the behaviour of dynamical processes with \textit{perturbations}.
We prove that a globally optimal state is stochastically stable under the logit response dynamics.
We also show that the better/best response dynamics converges to locally optimal states.
Numerical experiments confirm theoretical findings regarding some features of these dynamics.

The convergence and stability results indicate that a traffic state is led toward an efficient state by the evolutionary implementation scheme of the marginal cost pricing \citep{Sandholm2002,Sandholm2005,Sandholm2007} in the DSO game, in which toll level is adjusted according to the realised traffic state on a day-to-day basis.
Here, we further examine whether such an implementation scheme is essentially important for achieving efficient states.
To this end, we compare the scheme with another typical scheme, ‘fixed implementation scheme’, in which toll level is set optimally in advance according to a (known or target) optimal state.
Through the comparison, we found that under the evolutionary implementation scheme, \textcolor{black}{(i) the total travel time decreases smoother to an efficient traffic state, and (ii) a traffic state would reach a more efficient state.}
Numerical experiments are finally conducted to validate the theoretical findings.


The remainder of this paper is organised as follows.
In Section~\ref{Sec:DSOgame}, we define the DSO game.
Section~\ref{Sec:PropDSOgame} presents the convergence and stability properties in the DSO game under natural evolutionary dynamics.
Section~\ref{Sec:Solutions} shows solution algorithms based on the evolutionary dynamics; using these algorithms, we conduct numerical experiments in the DSO game.
Section~\ref{Sec:DUEFCP} compares the evolutionary implementation scheme with the fixed implementation scheme, theoretically and numerically. 
Section~\ref{Sec:Conclusion} concludes the paper.

\section{DSO game}\label{Sec:DSOgame}
In this section, we formulate a DSO traffic assignment problem that deals with atomic users as a strategic game, `DSO game'.
The game consists of atomic users travelling through the network (i.e. players), sets of available routes for the users (i.e. strategy sets), and dynamic loading models determining their travel times (i.e. utility functions).
After explaining these components, we define Nash equilibrium of the DSO game.

\subsection{Player, strategy and utility}
A general road network with many-to-many origin--destination (OD) pairs is herein considered.
The network consists of a set of nodes $\mathcal{N}$ and a set of directed links $\mathcal{L}$.
The sets of origin nodes and destination nodes are denoted by $\mathcal{N}_{o}$ and $\mathcal{N}_{d}$, respectively.
The set of all acyclic routes from node $a$ to node $b$ is denoted by $\mathcal{R}(a,b)$.
When these nodes are not connected, \textcolor{black}{the set becomes an empty set.}

The length of link $l\in\mathcal{L}$ is denoted by $L_{l}$.
The free-flow speed, backward wave speed and saturation flow rate are constant and denoted by $v_{l}$, $w_{l}$ and $q_{l}$, respectively.
\textcolor{black}{Each link has a bottleneck at the end of the link, and the bottleneck capacity of link $l$ is denoted by $\mu_{l}$ $(\leq q_{l})$.
This capacity represents the maximum possible outflow rate from the link.
This means that we deal with not only a homogeneous link where the capacity of every section is equal, but also a heterogeneous link where the capacities of entrance and exit sections are different}, i.e. \textcolor{black}{the shape of the fundamental diagram of each location within each link is triangular or trapezoidal.}
The values of saturation flow and bottleneck capacity of links can be set sufficiently large such that \textcolor{black}{vehicle queues do not occur on the links}\footnote{\textcolor{black}{At an intersection, this setting represents that the intersection just downstream of the links is undersaturated.}}; we refer to such links as `uncapacitated link', whereas we refer to the other links as `capacitated links'.
If necessary, we distinguish these two types of links.

Each atomic user is a player in the DSO game.
The set of users is denoted by $\mathcal{P}$, and the number of users is denoted by $|\mathcal{P}|$.
The origin, destination and departure times of user $i\in\mathcal{P}$ are denoted by $o_{i}$, $d_{i}$ and $s_{i}$, respectively.
These are given exogenously.
All users departing from the same origin have different departure times.
\textcolor{black}{Each user selects an acyclic route between his/her origin and destination, i.e. each acyclic route between the origin and destination is a strategy of the user.}
The set of available routes for a user $i\in\mathcal{P}$ is denoted by $\mathcal{R}_{i}~(=\mathcal{R}(o_{i},d_{i}))$.
This set includes \textcolor{black}{a special strategy $\phi_{i}$}, which represents that user $i$ selects no route (i.e. user $i$ is not assigned to the network).
A set of strategies for all users (i.e. strategy profile) corresponds to a route choice pattern of all users, which is referred to as a `route profile' or just a (traffic) state.
A route profile is denoted by a vector $\mathbf{r}\equiv \{ r_{1},\ldots,r_{i},\ldots,r_{|\mathcal{P}|} \}\in\mathcal{R}$ where $\mathcal{R}\equiv \mathcal{R}_{1}\times\cdots \times \mathcal{R}_{|\mathcal{P}|}$.
\textcolor{black}{Note that $\mathcal{R}$ is finite as we consider only acyclic routes and a given finite number of users.}
For any route profile $\mathbf{r}$, the route choices of the users other than user $i$ are denoted by $\mathbf{r}_{-i}\equiv \{ r_{1},\ldots,r_{i-1},r_{i+1},\ldots,r_{|\mathcal{P}|} \}$.
With this notation, we sometimes represent a profile $\mathbf{r}$ as $(r_{i},\mathbf{r}_{-i})$ to clearly state the route of user $i$.

We assume that the utility of each user is equal to the negative value of his/her route marginal social cost on the selected route. 
This is usually \textcolor{black}{interpreted as the situation} in which each user tries to minimise his/her generalised travel cost under the marginal cost or Pigouvian pricing scheme.
The marginal cost of a user who selects a route for a given route profile is defined as the sum of the route travel time of the user (private cost) and the change in the route travel times of the other users by assigning the user (external cost).
Mathematically, the utility of user $i$ who selects the route $r_{i}\in\mathcal{R}_{i}$ for a given route profile $\mathbf{r}_{-i}$ is represented by, 
\begin{align}
&U_{i}(r_{i}, \mathbf{r}_{-i}) = - C_{i}(r_{i}, \mathbf{r}_{-i}) - E_{i}(r_{i}, \mathbf{r}_{-i}),\label{Eq:Utility_Evolutionary_0}
\\
&\text{where}\quad E_{i}(r_{i}, \mathbf{r}_{-i}) = \sum_{i' \in \mathcal{P}\setminus \{i\}}\left\{ C_{i'}(r_{i}, \mathbf{r}_{-i}) - C_{i'}(\phi_{i}, \mathbf{r}_{-i})   \right\},\label{Eq:Utility_Evolutionary}
\end{align}
where $C_{i}(r_{i}, \mathbf{r}_{-i})$ and $E_{i}(r_{i}, \mathbf{r}_{-i})$ are private and external costs of user $i$, respectively, and $C_{i'}(\phi_{i}, \mathbf{r}_{-i})$ represents the travel time of user $i'$ when user $i$ is \textit{not} assigned to the network\footnote{The marginal costs can be evaluated exactly, unlike the conventional fluid approach wherein the evaluation is difficult and inexact~\citep[][]{Qian2012}; we compare the resulting travel times in the cases with the route profiles of $(r_{i}, \mathbf{r}_{-i})$ and $(\phi_{i}, \mathbf{r}_{-i})$.}.
We also denote by $TC(\mathbf{r})$ the total travel time of the transport system (referred to as total cost, hereinafter) for route profile $\mathbf{r}$, as follows: 
\color{black}
\begin{align}
TC(\mathbf{r}) = \sum_{i \in \mathcal{P}} C_{i}(\mathbf{r}),\quad \forall \mathbf{r}\in\mathcal{R}.\label{Eq:TTC}
\end{align}
\color{black}

\textcolor{black}{The travel time (and external cost) of each user for a given route profile should be uniquely determined according to a dynamic loading model.
To this end, we assume that the FIFO principle and causality~\citep[][]{Carey2003} are satisfied on each link.
On an intersection (including merge or diverge) node, we assume that a node model satisfies conditions consistent with generic requirements for macroscopic node models in \cite{Tampere2011} and \cite{Smits2015}, such as non-vehicle holding, an FIFO rule for diverging behaviour, to ensure natural user behaviour; we also assume that a priority rule for merging behaviour (i.e. which user should be prioritised when users on multiple incoming links can enter a link at the same time) is specified so that a user departure time pattern from the intersection is uniquely determined for a given user arrival time pattern.
Any appropriate model which satisfies the natural conditions for dynamic loading can be employed.}

\subsection{Nash equilibrium}
We employ the concept of pure Nash equilibrium as equilibrium of the users' route choices.
Mathematically, an equilibrium state $\mathbf{r}^{*}$ satisfies the following condition: 
\begin{align}
U_{i}(r^{*}_{i}, \mathbf{r}^{*}_{-i}) = \max_{r\in\mathcal{R}_{i}}U_{i}(r, \mathbf{r}^{*}_{-i}),\quad \forall i\in\mathcal{P}.\label{Eq:Nash-base}
\end{align}
The route of each user in a Nash equilibrium state is a best response to the route choices of all other users; we refer to the route as a `best response route'.
If each user has the unique best response route, $\mathbf{r}^{*}$ is called a \textit{strict} Nash equilibrium state.

\textcolor{black}{By substituting Eq.~\eqref{Eq:Utility_Evolutionary_0} into the condition~\eqref{Eq:Nash-base}, we obtain}: 
\begin{align}
\sum_{i' \in \mathcal{P}} C_{i'}(r^{*}_{i}, \mathbf{r}^{*}_{-i}) = \min_{r\in \mathcal{R}_{i}} \sum_{i' \in \mathcal{P}} C_{i'}(r, \mathbf{r}^{*}_{-i}),\quad \forall i\in\mathcal{P}.\label{Eq:Nash}
\end{align}

\noindent 
\textcolor{black}{This equation means that the total cost at equilibrium cannot be decreased by unilaterally changing routes.
Thus, each equilibrium state corresponds to a locally optimal state of a total cost minimisation problem\footnote{By definition, local maximum and saddle points do not become Nash equilibrium states, unlike in the fluid approach.}.
When the total cost of a locally optimal state is minimal among all locally optimal states, the optimal state is referred to as the globally optimal state.}

\textcolor{black}{Note that the globally optimal state always exists because the set of feasible states (i.e. route profiles) is finite; however, the uniqueness is not guaranteed in general since different route profiles could have the same total cost.}


\section{Convergence and stability in DSO games}\label{Sec:PropDSOgame}
This section analyses the convergence and stability in a DSO game.
We first show that a DSO game is a potential game.
By utilising this appealing property, we prove the convergence of evolutionary dynamics to locally optimal states.
Furthermore, we establish the stochastic stability of a globally optimal state under evolutionary dynamics with perturbations.

\subsection{DSO game is a potential game}
A potential game is a game wherein the change in the utility of a user, which results from a unilateral change in strategy, equals the change in the global utility referred to as a potential function.
This is formally defined as follows: 
\begin{defi}\textbf{(Potential game~\citep[][]{Monderer1996})}
A finite $n$-user game with action sets $\{ \mathcal{R}_{i}\}_{i=1}^{n}$ and utility functions $\{U_{i} \}_{i=1}^{n}$ is a potential game if and only if there is a function $\Pi$ : $\mathcal{R}\rightarrow \mathbb{R}$ such that
\begin{align}
U_{i}(r'_{i}, \mathbf{r}_{-i}) - U_{i}(r''_{i}, \mathbf{r}_{-i}) = \Pi(r'_{i}, \mathbf{r}_{-i}) - \Pi(r''_{i}, \mathbf{r}_{-i}), \quad \forall i\in\mathcal{P}, \forall \mathbf{r}_{-i}\in\mathcal{R}_{-i},\forall r',r''\in\mathcal{R}_{i}.
\end{align}
\end{defi}
\noindent According to this definition, we have the following theorem:

\begin{theo}
\textcolor{black}{A DSO game is a potential game whose potential function is the negative of the total cost function in Eq.~\eqref{Eq:TTC}}.
\end{theo}
\begin{prf}
Considering two route profiles $(r'_{i}, \mathbf{r}_{-i})$ and $(r_{i}, \mathbf{r}_{-i})$ where user $i$ changes his/her route from $r_{i}$ to $r'_{i}$, we have 
\begin{align*}
U_{i}(r'_{i}, \mathbf{r}_{-i}) - U_{i}(r_{i}, \mathbf{r}_{-i}) =& - C_{i}(r'_{i}, \mathbf{r}_{-i}) - E_{i}(r'_{i}, \mathbf{r}_{-i}) - \left\{- C_{i}(r_{i}, \mathbf{r}_{-i}) - E_{i}(r_{i}, \mathbf{r}_{-i})\right\} \notag \\
 = &- C_{i}(r'_{i}, \mathbf{r}_{-i}) + C_{i}(r_{i}, \mathbf{r}_{-i})- \sum_{i' \in \mathcal{P}\setminus \{i\}}\left\{ C_{i'}(r'_{i}, \mathbf{r}_{-i}) - C_{i'}(\phi_{i}, \mathbf{r}_{-i}) - C_{i'}(r_{i}, \mathbf{r}_{-i}) + C_{i'}(\phi_{i}, \mathbf{r}_{-i})     \right\}\notag \\
= &- \sum_{i\in\mathcal{P}} C_{i'}(r'_{i}, \mathbf{r}_{-i}) + \sum_{i\in\mathcal{P}} C_{i'}(r_{i}, \mathbf{r}_{-i}) = -TC(r'_{i}, \mathbf{r}_{-i}) + TC(r_{i}, \mathbf{r}_{-i}).
\end{align*}
Hence, the change in the utility of user $i$ equals the change in the negative of the total cost function.
This means that the negative of the total cost function is the potential function of this game; this proves the theorem.\qed
\end{prf}
\noindent It means that an improvement in the utility of each user always leads to an improvement in the \textcolor{black}{total cost}.

\subsection{Convergence}\label{Sec:DSOconvergence}
In the following, we consider a situation in which a DSO game is repeatedly played on a day-to-day basis.
In such a repeated game, at each time (i.e. day) $\tau \in \mathbb{N}$, each user $i\in\mathcal{P}$ takes the route $r^{\tau}_{i}\in\mathcal{R}_{i}$ and receives the utility $U_{i}(\mathbf{r}^{\tau})$ where $\mathbf{r}^{\tau}$ is the route profile at the time.
\textcolor{black}{On each day, one user who is going to change his/her route is selected randomly with equal probability $1/|\mathcal{P}|$, and the user calculates his/her utility of each route for a route profile on the previous day\footnote{\textcolor{black}{It would be difficult for the selected user to calculate external costs.
However, these costs are considered to be imposed by a road manager in the form of marginal cost or Pigouvian pricing schemes, as mentioned in Section~\ref{Sec:DSOgame}.
It follows that the user only have to calculate his/her private costs (i.e. route travel times) on the assumption that the other users keep taking their current routes, and this may not be so difficult.}}.
The other users must repeat their route choices from the previous day.
Then, based on the calculated utility, the selected user changes the route according to a behavioural rule common to all users by comparing the utilities of the current route and other routes. 
The resulting (day-to-day) evolution of the traffic state is called evolutionary dynamics.}


Let us investigate the convergence of the following two general evolutionary dynamics: better response and best response dynamics.
Under the better response dynamics, the selected user $i$ changes from route $r^{\tau}_{i}$ to $r^{\tau+1}_{i}$ if this strictly improves the user's utility compared with the previous day.
\textcolor{black}{If the user has multiple better response routes, the user selects one among them according to an arbitrary probability distribution\footnote{\textcolor{black}{For example, each user selects one better response route among multiple better response routes with equal probability or according to its utility.}}.}
If the user does not have a route satisfying this condition, the user does not change the route from that of the previous day.
Mathematically, the probability $p_{i}^{better}(r ;\mathbf{r}^{\tau})$ that selected user $i$ chooses route $r$ for a given route profile $\mathbf{r}^{\tau}$ is given as follows: 
\color{black}
\begin{align}
&p_{i}^{better}(r ;\mathbf{r}^{\tau})= 
\begin{cases}
1\quad & 
\text{if}\quad D_{i}(\mathbf{r}^{\tau}) = \emptyset \text{ \textcolor{black}{and} } r = r^{\tau}_{i},\\
q_{i}(r; \mathbf{r}^{\tau}) \quad
& \text{if}\quad D_{i}(\mathbf{r}^{\tau}) \neq \emptyset \text{ \textcolor{black}{and} } r\in D_{i}(\mathbf{r}^{\tau}),\\
0 & \text{otherwise},
\end{cases}\label{Eq:Prob_BetterResponse}
\end{align}
where $q_{i}(r; \mathbf{r}^{\tau})$ is the positive probability that the user $i$ chooses a better response route $r$ for a given route profile $\mathbf{r}^{\tau}$. The probability should satisfy the normalised condition, $\sum_{r\in D_{i}(\mathbf{r}^{\tau})}q_{i}(r; \mathbf{r}^{\tau}) = 1$, where $D_{i}(\mathbf{r}^{\tau})$ is the set of better responses of the user $i$ for the given route profile $\mathbf{r}^{\tau}$ and defined as follows:
\begin{align}
D_{i}(\mathbf{r}^{\tau}):=\left\{ r^{*}_{i}\ | \ r^{*}_{i}\in\mathcal{R}_{i}~\text{s.t.}~U_{i}(r^{*}_{i},\mathbf{r}^{\tau}_{-i}) > U_{i}(r^{\tau}_{i},\mathbf{r}^{\tau}_{-i}) \right\}.\label{eq:Better-response}
\end{align}
\color{black}
Then, the transition probability from a route profile $\mathbf{r}^{\tau}$ to $\mathbf{r}^{\tau+1}$ is represented as follows: 
\begin{align}
&p^{0}_{\mathbf{r}^{\tau}\mathbf{r}^{\tau+1}} = 
\begin{cases}
\cfrac{1}{|\mathcal{P}|}\cdot p^{better}_{i}(r_{i}^{\tau+1}; \mathbf{r}^{\tau})   
&\quad \text{if}\quad r_{i}^{\tau+1}\neq r_{i}^{\tau}\text{ \textcolor{black}{and} }\mathbf{r}^{\tau+1} = (r_{i}^{\tau+1}, \mathbf{r}_{-i}^{\tau}), \quad \forall i\in\mathcal{P},\\
\sum_{i\in\mathcal{P}} \cfrac{1}{|\mathcal{P}|}\cdot p^{better}_{i}(r^{\tau+1}_{i} ; \mathbf{r}^{\tau}) 
&\quad \text{if}\quad \mathbf{r}^{\tau+1} = \mathbf{r}^{\tau},\\
0 &\quad \text{otherwise}.
\end{cases}\label{Eq:Better-Markov}
\end{align}
The first case shows that one user $i$ changes his/her route $r_{i}^{\tau}$ to a different route $r_{i}^{\tau + 1}$;
the second case shows that a user does not change the route;
the third case shows that two or more users change their routes simultaneously, and this probability must be zero.
The transition probability suggests that the probability distribution of route profiles at $\tau+1$ depends only on the route profile at $\tau$.
Thus, the resulting stochastic process of route profiles $\{\mathbf{r}^{\tau} \}_{\tau\in \mathbb{N}}$ becomes a Markov chain that describes the transition probability between each route profile, and its state space is the set of route profiles $\mathcal{R}$.


As a direct consequence of the fact that a DSO game is a potential game, we can establish the following theorem regarding the convergence of the dynamics: 
\begin{theo}\label{Theo:ConBetter}{\rm \textbf{(Convergence of better response dynamics).}
In a DSO game under the better response dynamics, a route profile converges almost surely to a Nash equilibrium state (locally or globally optimal state) from an arbitrary initial profile, and this state has a \textcolor{black}{lower total cost} than that of the initial state.}
\end{theo}
\begin{prf}
\textcolor{black}{Since the game is a potential game, the better response dynamics which increases the potential of the DSO game converges to a state maximising the potential locally or globally as $\tau\rightarrow \infty$}\footnote{\textcolor{black}{A traffic state might get stuck into non-Nash equilibrium if a user that cannot change his/her route is continuously given opportunity to change the route and the potential does not increase.
However, it is obvious that such a probability converges to zero as $\tau\rightarrow \infty$.}}.
Moreover, since the potential is equal to the negative value of the \textcolor{black}{total cost}, an increase in the potential by the better responses always decreases the \textcolor{black}{total cost}.
Thus, the proposition is proved.\qed
\end{prf}

Under the best response dynamics, the selected user changes his/her route to maximise the utility. 
In other words, given any current route profile $\mathbf{r}^{\tau}$, the user $i$ chooses a route randomly from the following route set $B_{i}(\mathbf{r}^{\tau})$: 
\begin{align}
B_{i}(\mathbf{r}^{\tau}) := \left\{ r^{*}_{i}\ |\ r^{*}_{i} \in\mathcal{R}_{i}~\text{s.t.}~ 
\max_{r\in\mathcal{R}_{i}}U_{i}(r,\mathbf{r}^{\tau}_{-i})   \right\}.
\end{align}
Then, the probability $p_{i}^{best}(r ;\mathbf{r}^{\tau})$ that selected user $i$ chooses route $r$ for a given route profile $\mathbf{r}^{\tau}$ is given as follows:
\begin{align}
&p_{i}^{best}(r ;\mathbf{r}^{\tau})= 
\begin{cases}
\cfrac{1}{|B_{i}(\mathbf{r}^{\tau})|} \quad 
& \text{if}\quad r\in B_{i}(\mathbf{r}^{\tau}),\\
0 & \text{otherwise}.
\end{cases}\label{Eq:Prob_BestResponse}
\end{align}
Transition probability between route profiles can be derived by replacing the route choice probability $p^{better}_{i}(r^{\tau+1};\mathbf{r}^{\tau})$ in Eq.~\eqref{Eq:Better-Markov} by $p^{best}_{i}(r^{\tau+1};\mathbf{r}^{\tau})$.
Thus, the stochastic process $\{\mathbf{r}^{\tau}\}_{\tau\in \mathbb{N}}$ under the best response dynamics is also a Markov chain with finite state space $\mathcal{R}$.

The best response dynamics has the following property compared to the better response dynamics: 
a user taking his/her best response route can change the route to another best response route if it exists, i.e. users can switch their routes among their best response routes\footnote{\textcolor{black}{The better response dynamics does not allow switching among the best response routes because users cannot change their routes to other routes whose utility is the same.}}.
This means that even when the current route profile is Nash equilibrium, a route profile can deviate from that equilibrium if some user finds another route whose marginal cost is the same with that of the current route.
This property yields the following corollary about a closed communication class, which is an absorbing set of states in a Markov chain, of the best response dynamics\footnote{\textcolor{black}{A communication class $\mathcal{C}\subseteq \mathcal{R}$ is a set of states (route profiles) whose members can reach each other (i.e. communicate) under the corresponding evolutionary dynamics, and no state in $\mathcal{C}$ communicates with any state outside $\mathcal{C}$.
Also, a communication class $\mathcal{C}$ is said to be closed if no state in $\mathcal{C}$ reaches any state outside $\mathcal{C}$.
For details, see for example, \cite{Meyn2009}.}}: 


\begin{coro}\label{Coro:ConBest}
In a DSO game, any closed communication class in a Markov chain generated by the best response dynamics consists of Nash equilibrium states with equal \textcolor{black}{total costs}.
\end{coro}
\begin{prf}
If two traffic states have different \textcolor{black}{total costs}, these states are not reachable from each other by the best response dynamics; the state with the higher total cost is not reachable from the other.
It follows that all the traffic states in any communication class should have the same total cost, in order to be able to `communicate'.

Moreover, if there exists a traffic state that is not a Nash equilibrium state (i.e. there exists at least one user who could improve the utility) in a communication class, it is not a closed class: if the \textcolor{black}{total cost} is improved by a best response of a user and a traffic state deviates from the communication class, the best response dynamics would never return to traffic states in that class.
It follows that all closed communication classes consist of Nash equilibrium states with equal \textcolor{black}{total costs}.
Thus, the corollary is proved. \qed
\end{prf}

\noindent \textcolor{black}{This corollary means that a route profile could change among a set of states with equal total costs.}
From this corollary, we obtain the following theorem showing the convergence of the best response dynamics: 
\begin{theo}{\rm \textbf{(Convergence of best response dynamics).}}
\textcolor{black}{In a DSO game under the best response dynamics, a route profile converges almost surely to \textit{a set of} Nash equilibrium states with equal total costs from an arbitrary initial profile, and these states have a total cost lower than that of the initial state.}
\end{theo}
\begin{prf}
All Markov chains eventually \textcolor{black}{enter} a closed communication class with probability 1 as $\tau\rightarrow \infty$.
From this fact and \textbf{Corollary~\ref{Coro:ConBest}}, it follows that a route profile generated by the best response dynamics is absorbed into a set of Nash equilibrium states with equal \textcolor{black}{total costs}.
Moreover, as in the statement in the proof of \textbf{Theorem~\ref{Theo:ConBetter}}, the \textcolor{black}{total cost} is lower than that of the initial state.\qed

\end{prf}

It might look undesirable from a convergence perspective that a traffic state could constantly change under the best response dynamics.
However, such behaviour of the dynamics yields the following good property: a route profile can deviate from an inefficient (or locally) optimal state where the better response dynamics could converge.
\textcolor{black}{Specifically, under the best response dynamics, a route profile can shift from an equilibrium state to another state with an equal \textcolor{black}{total cost} since users can switch their routes among their best response routes, as mentioned above.
Such a switch could cause the following ripple effect~\citep[][]{Satsukawa2019}: route change of a user affects utilities of other users and makes their routes into non-best response routes.
This means that these users can find routes whose marginal costs are lower than those they currently select after the shift.
Thus, if such users are given opportunities to change their routes, the \textcolor{black}{total cost} can decrease by their best responses.}

Meanwhile, under the better response dynamics, a route profile cannot deviate from any equilibrium state because a user can only select a route whose utility is \textit{strictly} higher.
In other words, in a Markov chain generated by the better response dynamics, all Nash equilibrium states are regarded as absorbing states of the Markov chain.
\textcolor{black}{It is thus expected that total costs generated by the best response dynamics tend to be lower than those by the better response dynamics.
In other words, relaxing such strict improvement property in user utility would improve the efficiency of a transport system.}


\begin{figure}[t]
	\begin{center}
	\hspace{0mm}
    \includegraphics[width=0.9\linewidth,clip]{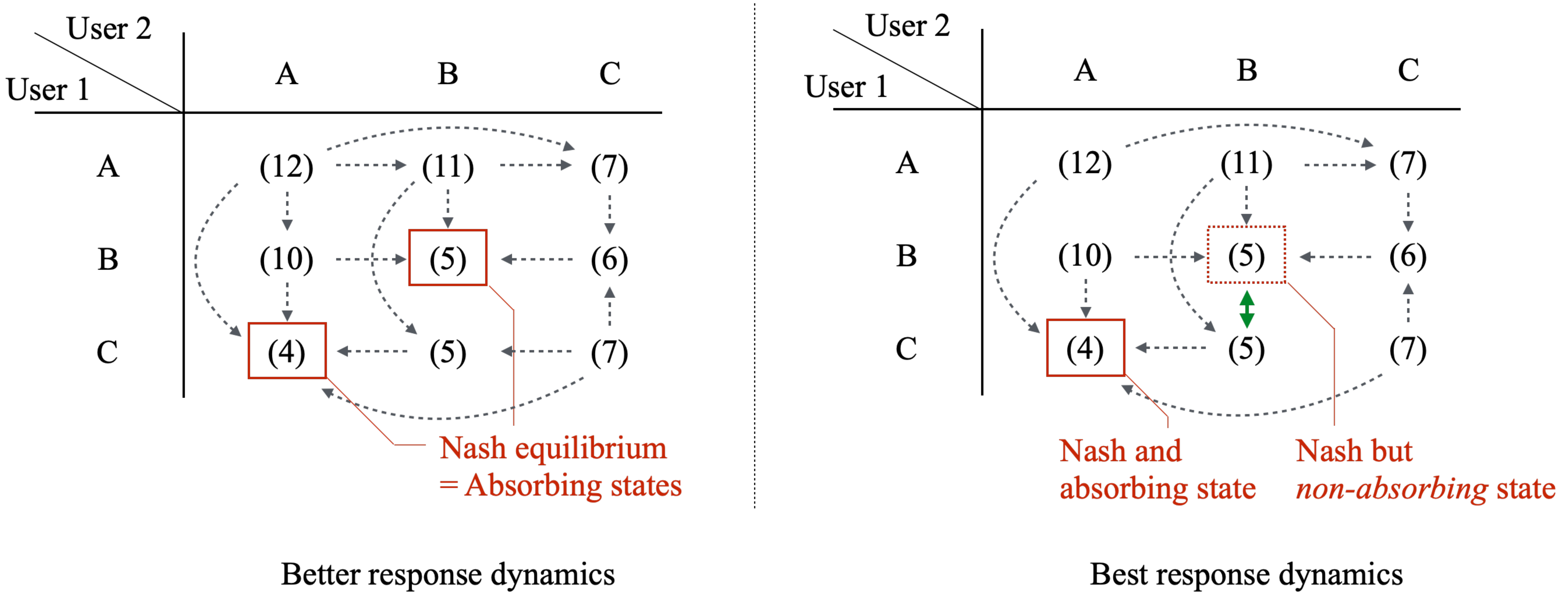}
	\end{center}
    \vspace{-4mm}
	\caption{\textcolor{black}{Matrices of a two-user DSO game in which each number is the \textcolor{black}{total cost} corresponding to each route profile. Each dotted arrow represents a better or best response.}}
    \vspace{-3mm}
    \label{Fig:PayoffMat}
\end{figure}

An example demonstrating the difference is illustrated in Figure~\ref{Fig:PayoffMat}.
Each matrix in the figure shows the total costs of the route profiles and the better or best responses in a two-user DSO game.
There exist two Nash equilibrium states $(B, B)$ and \textcolor{black}{$(C, A)$}\footnote{We can identify which route profiles are Nash equilibrium states from the matrix showing the \textcolor{black}{total cost}: since a DSO game is a potential game whose potential function corresponds to the \textcolor{black}{total cost}, improvement in the \textcolor{black}{total cost} by a route change of a user always leads to an improvement in the utility of the user.}.
As shown in the left matrix, if the better response dynamics is employed, the route profile can converge to either of the equilibrium states since all equilibrium states are absorbing states of the dynamics.
In contrast, as shown in the right matrix, if the best response dynamics is employed, a route profile can deviate from $(B, B)$ and change to \textcolor{black}{$(C, B)$} since user 1 need not strictly improve the utility.
As a result, the route profile converges to \textcolor{black}{$(C, A)$}, whose \textcolor{black}{total cost} is lower than that of \textcolor{black}{$(C, B)$.}


\textcolor{black}{Although we here show the convergence of the deterministic evolutionary dynamics, it is essentially difficult to establish the convergence of such dynamics to the globally optimal state, as mentioned in Section~\ref{Sec:Introduction}.
In order to avoid getting stuck in (inefficient) locally optimal states where the best response dynamics also converges, it is necessary to introduce \textit{perturbations} that lead traffic states to more efficient states by allowing the shift to worse traffic states with higher total costs.
The next section examines evolutionary dynamics with perturbations, and shows the \textit{stochastic stability of globally optimal states}, which is regarded as the convergence concept of such perturbed dynamics.
}

\subsection{Stability}

\subsubsection{Definition of stochastic stability}
The following is a brief summary of \cite{Young1993}.
We first consider a finite state Markov chain over the state space $\mathcal{R}$ generated by evolutionary dynamics without perturbation (e.g. best response dynamics); we refer to this kind of evolutionary dynamics and the Markov chain as `unperturbed evolutionary dynamics' and `unperturbed Markov chain', respectively.
Let $\mathbf{P}^{0}$ be the transition matrix of the Markov chain.
Further, let also $p^{0}_{\mathbf{r}\mathbf{r}'}$ be an element of the matrix which represents the transition probability from state (i.e. route profile) $\mathbf{r}\in\mathcal{R}$ to state $\mathbf{r}'\in\mathcal{R}$.

We then consider another Markov chain that is generated by a perturbed version of the (unperturbed) evolutionary dynamics under which users are subjected to a small perturbation whose size is indexed by a scalar $\epsilon$, referred to as a `perturbed Markov chain'.
$\epsilon$ takes on all values in some interval $(0,a]$.
Let $\mathbf{P}^{\epsilon}$ be the corresponding transition probability matrix.
We assume that the perturbed Markov chain satisfies the following conditions:
\begin{align}
&\mathbf{P}^{\epsilon} \text{ is aperiodic and irreducible (i.e. the perturbed Markov chain is ergodic) for all $\epsilon\in(0,a]$,}\label{Eq:RPMC-0}\\
&\lim_{\epsilon\rightarrow 0^{+}} p^{\epsilon}_{\mathbf{r}\mathbf{r}'} = p^{0}_{\mathbf{r}\mathbf{r}'}\label{Eq:RPMC-1}\\
&\text{and } p^{\epsilon}_{\mathbf{r}\mathbf{r}'} > 0 \text{ for some $\epsilon$ implies } 
\exists c(\mathbf{r}\rightarrow \mathbf{r}')\geq 0 
\text{ s.t. } 0 <  \lim_{\epsilon\rightarrow 0^{+}}\epsilon^{-c(\mathbf{r}\rightarrow \mathbf{r}')}p^{\epsilon}_{\mathbf{r}\mathbf{r}'} < \infty.\label{Eq:RPMC-2}
\end{align}
Condition~\eqref{Eq:RPMC-0} implies that the perturbed Markov chain has a unique stationary distribution $\boldsymbol{\pi}^{\epsilon}$ for every $\epsilon$\textcolor{black}{~\citep[pp.241-243,][]{Stewart2009}}\footnote{\textcolor{black}{These properties are violated when, for example, there exists a zero probability path between a pair of some states, say $i$ and $j$, i.e. state $i$ is not reachable from state $j$.}}.
Condition~\eqref{Eq:RPMC-1} implies that $\mathbf{P}^{\epsilon}$ converges to the unperturbed one, $\mathbf{P}^{0}$.
Condition~\eqref{Eq:RPMC-2} implies that $\mathbf{P}^{\epsilon}$ approaches $\mathbf{P}^{0}$ at an exponentially smooth rate. 
The scalar $c(\mathbf{r}\rightarrow \mathbf{r}')$ is called the `resistance' of transition $\mathbf{r}\rightarrow \mathbf{r}'$. 
This value represents the degree of intensity of mistakes required for this transition (e.g. the minimum number of mistakes). 
If the transition is allowed under the unperturbed evolutionary dynamics, $c(\mathbf{r}\rightarrow \mathbf{r}') = 0$.
The specific Markov chain satisfying \eqref{Eq:RPMC-0}-\eqref{Eq:RPMC-2} is called a `regular perturbed Markov chain' of $\mathbf{P}^{0}$.

The unperturbed Markov chain may have multiple stationary distributions, whereas the (regular) perturbed Markov chain has a unique stationary distribution.
Such a unique stationary distribution converges to one of the stationary distributions of the unperturbed Markov chain, as $\epsilon\rightarrow 0$.
This means that the perturbations effectively select one stationary distribution.
This stationary distribution yields the observation probability of each state when the process of the Markov chain with the perturbations runs for a long time.
Stochastic stability is then defined as follows: 
\begin{defi}\textbf{(Stochastic stability~\citep[][]{Young1993}).}
A state $\mathbf{r}\in\mathcal{R}$ is `stochastically stable' relative to the Markov chain $\mathbf{P}^{\epsilon}$ if \textcolor{black}{$\lim_{\epsilon\rightarrow 0} \pi^{\epsilon}_{\mathbf{r}} > 0$.}
\end{defi}
\noindent Over the long run, states that are not stochastically stable are observed infrequently compared to states that are stochastically stable, provided that the probability of mistakes is small as $\epsilon\rightarrow 0$.

\subsubsection{Stochastic stability of globally optimal states}
We establish the stability of global optimal states by employing \textit{logit response dynamics}~\citep[][]{Blume1993}.
Under the dynamics, \textcolor{black}{the probability $p_{i}^{\beta}(r,\mathbf{r}^{\tau})$ that a selected user $i$ chooses route $r$ for a given route profile $\mathbf{r}^{\tau}$ is given as follows:}
\begin{align}
p_{i}^{\beta}(r,\mathbf{r}^{\tau}) = \cfrac{\exp(\beta U_{i}(r,\mathbf{r}_{-i}))}{\sum_{r'\in\mathcal{R}_{i}}\exp(\beta U_{i}(r',\mathbf{r}_{-i}))}\label{Eq:LogitDynamics}
\end{align}
where $\beta \in (0,\infty)$ measures the degree of noise in the best response.
Note that the logit response dynamics is a perturbed best response dynamics, and converges to the best response dynamics when $\beta\rightarrow \infty$ (this corresponds to $\epsilon\rightarrow 0$).

We then show that a globally optimal state is a stochastically stable under the logit response dynamics:

\begin{theo}{\rm \textbf{(Stochastic stability of logit response dynamics).}\label{Theo:SSS-DSO}
Consider the logit response dynamics in a DSO game.
The stochastically stable state is the globally optimal state minimising the total cost.}
\end{theo}
\begin{prf}
\cite{Marden2012b} show that the stochastically stable states in a potential game with the logit response dynamics are the set of potential maximisers.\qed
\end{prf}
\noindent \textcolor{black}{This theorem means that, with the help of the stochastic perturbation, a traffic state avoids converging to a locally (and inefficient) optimal states, and approaches the state minimising the \textcolor{black}{total cost} over the long run.}


\section{Solution algorithms and numerical experiments}\label{Sec:Solutions}
In this section, we briefly summarise the evolutionary dynamics with and without perturbations as solution algorithms for the DSO game.
We then show numerical experiments that demonstrate the theoretical properties of the evolutionary dynamics shown in the previous section.

\subsection{Solution algorithms}
We first show deterministic algorithms for computing a \textcolor{black}{\textit{locally optimal state}} based on the better and best response dynamics, as follows. 

\vspace{2mm}
\begin{algorithm}[H]\label{Algo:DeterministicSAlgo}
\caption{\textbf{Deterministic solution algorithms with better or best response dynamics}}

{
0. \textbf{Initialisation}: \textcolor{black}{Select either of the better or best response dynamics.}
Set $\tau = 0$ where $\tau$ represents the iteration counter.
Set also an initial route profile $\mathbf{r}^{0}$.

1. \textbf{Select a user conducting a better/best response}: 
Randomly select one user $i\in\mathcal{P}$ who can conduct a better or best response.
Calculate the utility of the user for each route from the dynamic loading model.
\textcolor{black}{Change the route of the selected user to a new route $r^{*}_{i}\in\mathcal{R}_{i}$ according to the selected dynamics (Eq.~\eqref{Eq:Prob_BetterResponse} or Eq.~\eqref{Eq:Prob_BestResponse}).}

2. \textbf{Update the route profile and judge the convergence}: Update the route profile: $\mathbf{r}^{\tau+1} = (r^{*}_{i},\mathbf{r}_{-i})$. 
If $\mathbf{r}^{\tau+1}$ is not a Nash equilibrium state, let $\tau:= \tau+1$ and go back to Step 1.
If $\mathbf{r}^{\tau+1}$ is a Nash equilibrium state, then terminate the algorithm.
}
\end{algorithm}
\vspace{2mm}

Next, we show a stochastic algorithm for computing a \textcolor{black}{\textit{globally optimal state}}, based on the logit response dynamics.
It is clear that as $\beta \rightarrow \infty$, the stationary distribution of the logit response dynamics in a DSO game converges to states maximising the potential function.
However, this does not imply the convergence of the logit response dynamics if $\beta$ is chosen as fixed.
To derive a globally optimal state, we have to develop an algorithm based on the logit resonse dynamics with a time-dependent perturbation parameter $\beta (\tau)$ that guarantees the convergence of the resulting stationary distribution to the potential maximiser, i.e. the state minimising \textcolor{black}{total cost}.
The framework of the stochastic algorithm is described as follows:

\vspace{2mm}
\begin{algorithm}[H]\label{Algo:StochasticSAlgo}
\caption{\textbf{Stochastic solution algorithm with logit response dynamics}}

{
0. \textbf{Initialisation}: Set $\tau = 0$ where $\tau$ represents the iteration counter.
Set also an initial route profile $\mathbf{r}^{0}$.

1. \textbf{Select a user changing the route according to the logit response dynamics}: 
\textcolor{black}{Randomly select one user $i\in\mathcal{P}$ and calculate the utility of the user for each route.
Change the user's route to a new route $r \in \mathcal{R}_{i}$ according to Eq.~\eqref{Eq:LogitDynamics}.}

2. \textbf{Update the route profile and perturbation parameter}: Update the route profile: $\mathbf{r}^{\tau+1} = (r,\mathbf{r}_{-i})$. 
Update also the perturbation parameter according to a specified schedule: $\beta: = \beta(\tau)$.

3. Go back to Step 1 and repeat.
}
\end{algorithm}
\vspace{2mm}

\noindent Regarding the condition of the perturbation parameter in the algorithm guaranteeing the convergence to a globally optimal state, we show a proposition based on the theorem presented in \cite{Tatarenko2014}, as follows:

\begin{prop}\label{Theo:LogitPara}
Consider a DSO game with $|\mathcal{P}|$ users.
Then, \textbf{Algorithm~\ref{Algo:StochasticSAlgo}} with $\beta(\tau) = \ln(\tau + 1)/(\lceil |\mathcal{P}|/2 \rceil)$ guarantees the probabilistic convergence of route profiles of the DSO game to the maximisers of the potential function, i.e.:
\begin{align}
\lim_{\tau\rightarrow \infty} \Pr \{   \mathbf{r}^{\tau} \in \{  \mathbf{r}^{*}  \mid     \Pi(\mathbf{r}^{*}) = \max_{\mathbf{r}} \Pi(\mathbf{r})   \}  \} = 1.
\end{align}
\end{prop}

\begin{prf}
See \ref{App:LogitPara}. \qed
\end{prf}

\noindent Although this proposition provides a setting of the perturbation parameter that guarantees convergence to a globally optimal state, the perturbation parameter has a logarithmic dependence on time. 
This means that the convergence speed is very low.
Thus, in the next subsection addressing numerical experiments of DSO games, we also introduce a more practical setting, in which $\beta(\tau)$ is a polynomial function, and test its performance.

\begin{figure}[t]
	\begin{center}
	\hspace{0mm}
    \includegraphics[width=90mm,clip]{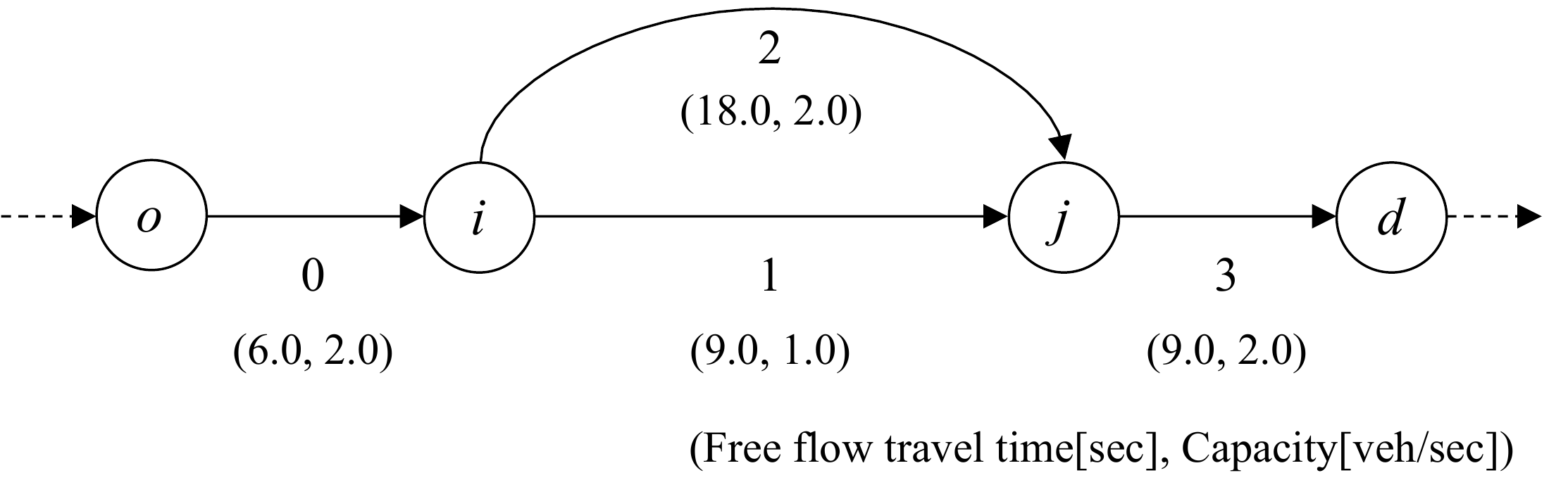}
	\end{center}
    \vspace{-4mm}
	\caption{Simple network with two parallel routes}
    \vspace{-0mm}
    \label{Fig:Network}
\end{figure}

\begin{figure}
	\begin{minipage}[t]{0.49\textwidth}
		\centering
		\includegraphics[clip, width=0.8\columnwidth]{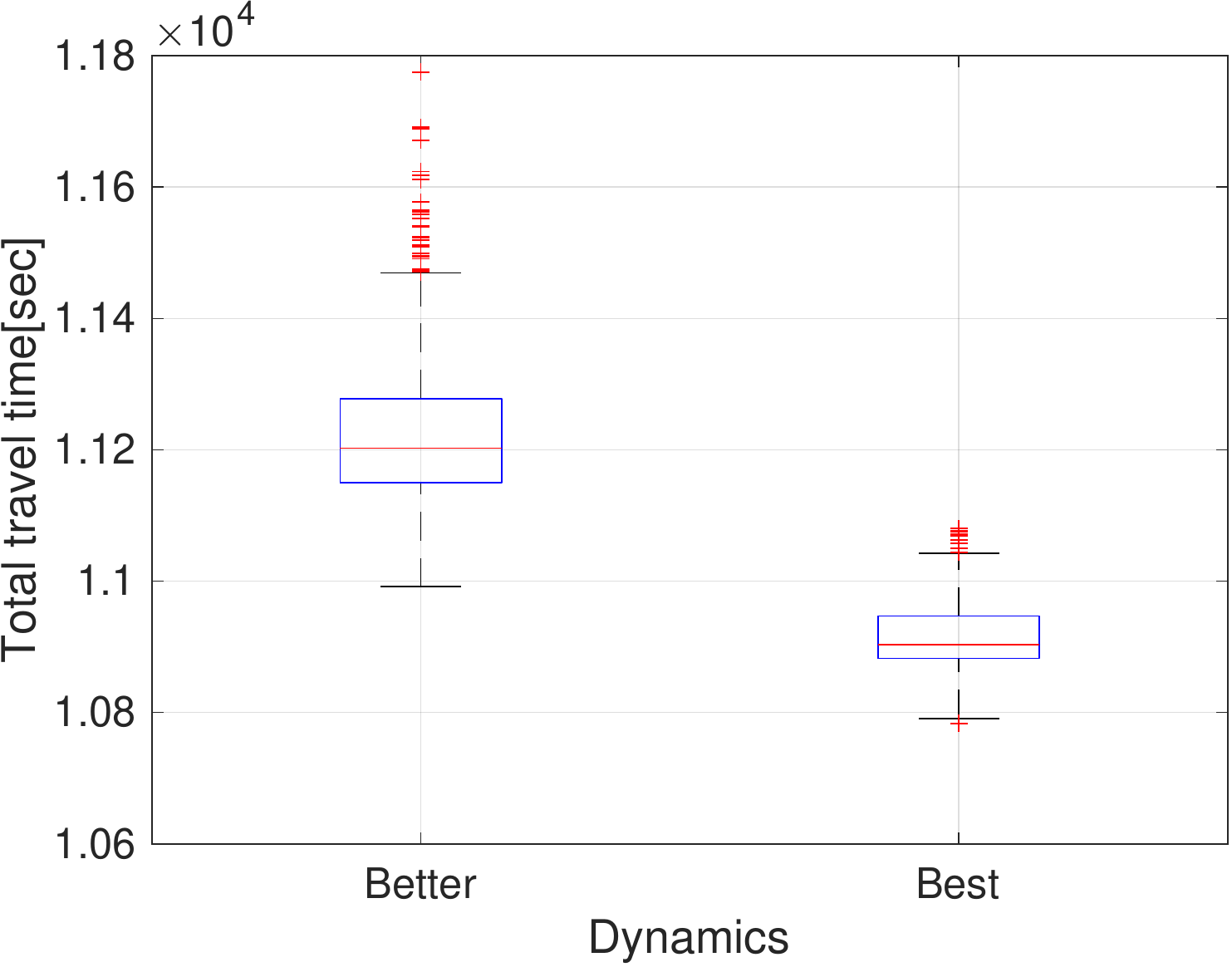}
		\subcaption{\textcolor{black}{Total costs} of the best optimal states in sample paths}\label{Fig:TTTDSO_D}
	\end{minipage}
	\begin{minipage}[t]{0.49\textwidth}
		\centering
		\includegraphics[clip, width=0.9\columnwidth]{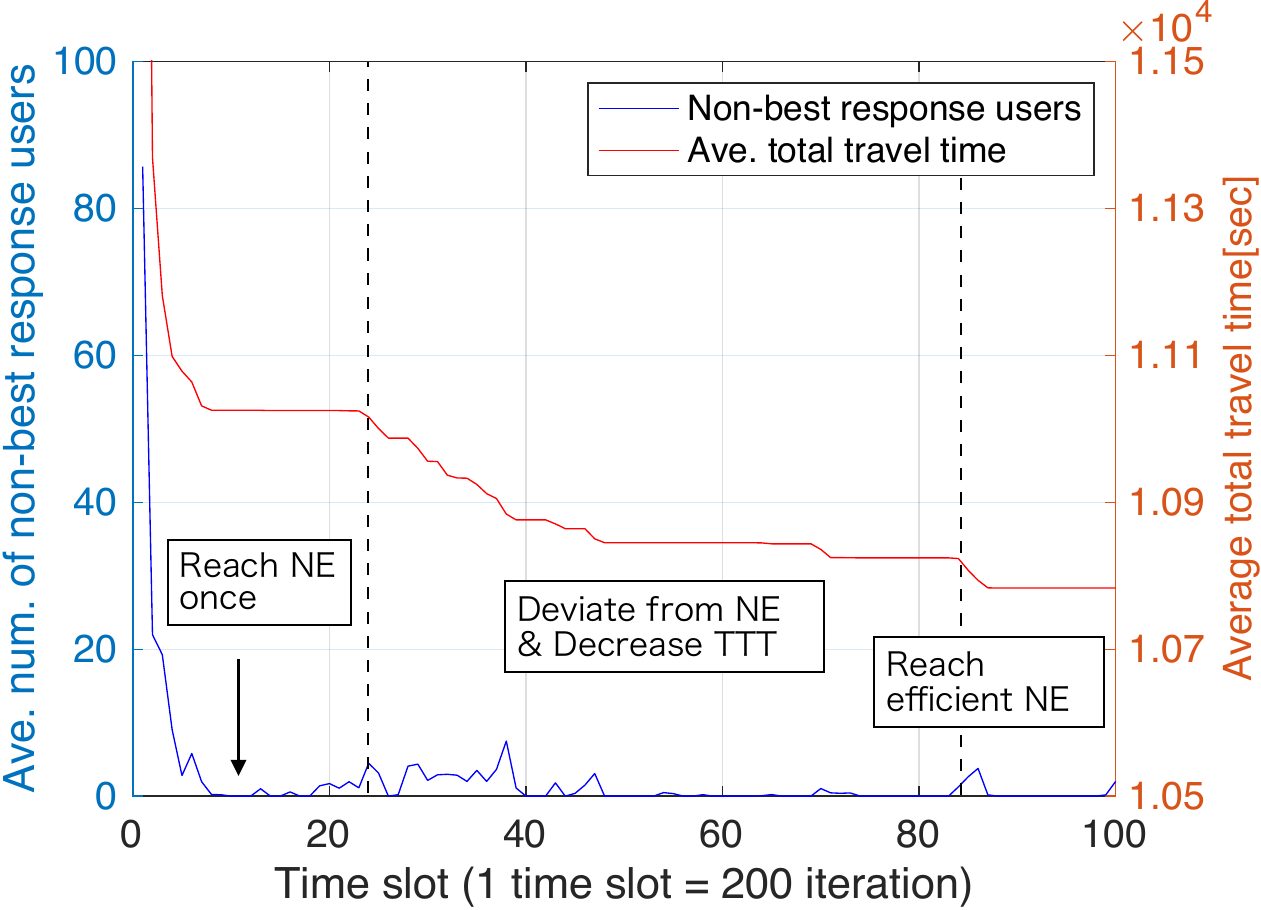}
		\subcaption{Convergence process of best response dynamics}\label{Fig:DSO_BestCon}
	\end{minipage}
	\vspace{-1mm} 
	\caption{Results regarding deterministic solution algorithms}
	\label{Fig:Results_D_simple}
	\vspace{-0mm}
\end{figure}

\begin{figure}[t]
	\begin{center}
	\hspace{0mm}
    \includegraphics[width=70mm,clip]{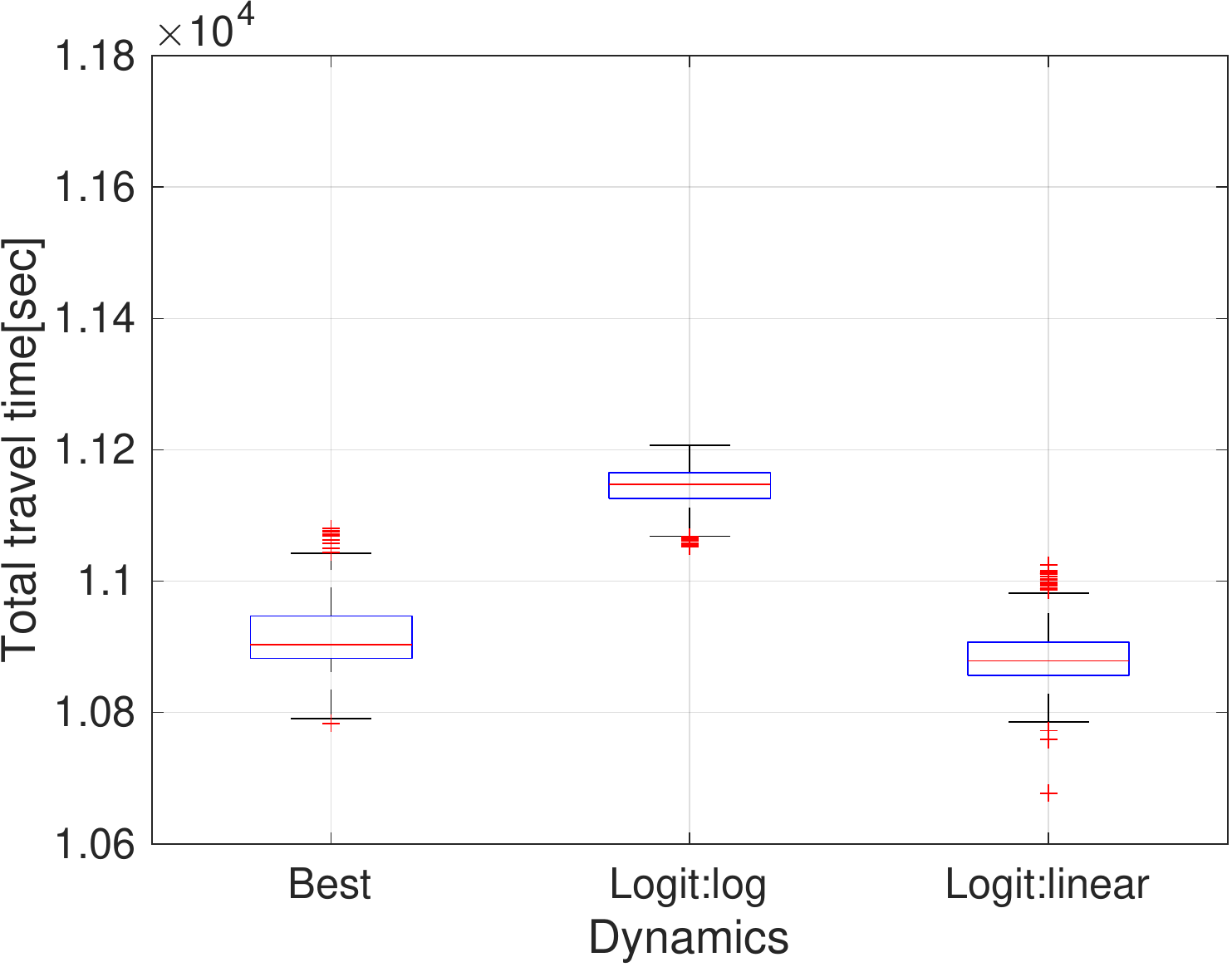}
	\end{center}
    \vspace{-4mm}
	\caption{Minimum \textcolor{black}{total costs} obtained from each solution algorithm}
    \vspace{-0mm}
    \label{Fig:TTTDSO_S}
\end{figure}

\subsection{\textcolor{black}{Numerical experiments}}

\subsubsection{Simple network}
We first consider a simple network with a single OD pair, as shown in Figure~\ref{Fig:Network}.
There exists two parallel routes (links) and each link has a bottleneck section with a bottleneck capacity at the end.
The physical conditions of each link (e.g. free-flow travel time and capacity) are summarised in the figure.
\textcolor{black}{Note that the saturation flow rates of all links are set to $\SI{2.0}{veh/sec}$.}
Route 1 is the shortest-distance route and consists of links $0$, $1$ and $3$, whereas route 2 is referred to as the bypass route and consists of links $0$, $2$ and $3$.
The total number of users is 400, and they depart from the origin with fixed time headway, $\SI{0.5}{sec/veh}$.

\textcolor{black}{For the numerical experiments, it is necessary to specify a dynamic loading model for calculating the user utility.
In this study, we consider a dynamic loading model which consists of two sub models: a link model and a node model.
A link model determines two kinds of `possible' times on each link according to existing vehicle trajectories and a specified car-following model: earliest possible departure times that vehicles can depart from the link when the downstream links are empty; earliest arrival times that vehicles can enter the link while satisfying the time-headway restriction required for the car-following behaviour.
In other words, the link model determines boundary conditions for determining departure and arrival times of vehicles on each link.
A node model then calculates `actual' departure and arrival times of vehicles on links such that the actual departure times are consistent with the boundary conditions of upstream and downstream links of each node.
The framework of the dynamic loading model is consistent with the `demand/supply approach' proposed by \cite{Daganzo1994,Daganzo1995} (see \ref{App:DNL} for the details of link and node models).}

In the numerical experiment, we compare the \textcolor{black}{total costs} derived by the solution algorithms (with better, best and logit responses) to observe the differences in the convergence properties between them.
For each algorithm, we generate 1,000 sample paths with different initial route profiles; the number of iterations for each sample path is set as 20,000 times, which is enough for each algorithm to find at least a locally optimal state (i.e. equilibrium).
\textcolor{black}{For the better response dynamics, we assume that if there exist multiple better response routes, each user selects one better response route among them with equal probability, i.e. $q_{i}(r;\mathbf{r}^{\tau}) = 1/|D_{i}(\mathbf{r}^{\tau})|$ in Eq.~\eqref{Eq:Prob_BetterResponse}.}
Note that we consider two different time-dependent perturbation parameters, logarithmic and linear decreasing, for the logit response dynamics. 
Specifically, in the former, the perturbation parameter at iteration $\tau$ is set as \textcolor{black}{$\ln(\tau+1)/200$}; in the latter, the parameter is set as $(\tau+1)/5000$.

Figure~\ref{Fig:TTTDSO_D} shows the distribution of the \textcolor{black}{total cost} of the best optimal state obtained in each sample path by the deterministic solution algorithms. 
From this figure, we see that there are multiple optimal (i.e. equilibrium) states with different \textcolor{black}{total costs} in this DSO game.
The difference between the worst and the best cases is approximately \SI{9}{\%}.

\textcolor{black}{This figure also shows that the total costs generated by the best response dynamics tend to be lower than those by the better response dynamics.
This is because, as mentioned in Section~\ref{Sec:DSOconvergence}, the behaviour of the two evolutionary dynamics is different near equilibrium even in this simple network.
Specifically, when there exist users whose utilities of the two routes are the same (i.e. both of the routes are the best response routes), the users can switch their routes among the two best response routes under the best response dynamics.
Then, there is a possibility that a route profile under the best response dynamics shifts from Nash equilibrium to more efficient equilibrium.}

Let us look at Figure~\ref{Fig:DSO_BestCon} that shows the behaviour of the best response dynamics in more detail.
In this figure, for every time slot (200 iterations), we calculate the average \textcolor{black}{total costs}, and number of users who do not take their best response routes (referred to as `non-best response users') in the route profile.
\textcolor{black}{From this figure, we see that the number of non-best response users increases from zero during the process; then, the \textcolor{black}{total cost} sometimes decreases in response to the increase.
This result comes from ripple effects caused by best responses.
Specifically, a user's best response among the two (best response) routes affects utilities of other users and makes their route into non-best response routes; then, if one of such users takes the best response, the traffic state moves to a more efficient state with deviating from an inefficient one. 
This result thus confirms the properties of the dynamics that the theory predicts.}



\begin{figure}
	\begin{minipage}[t]{0.49\textwidth}
		\centering
		\includegraphics[clip, width=0.9\columnwidth]{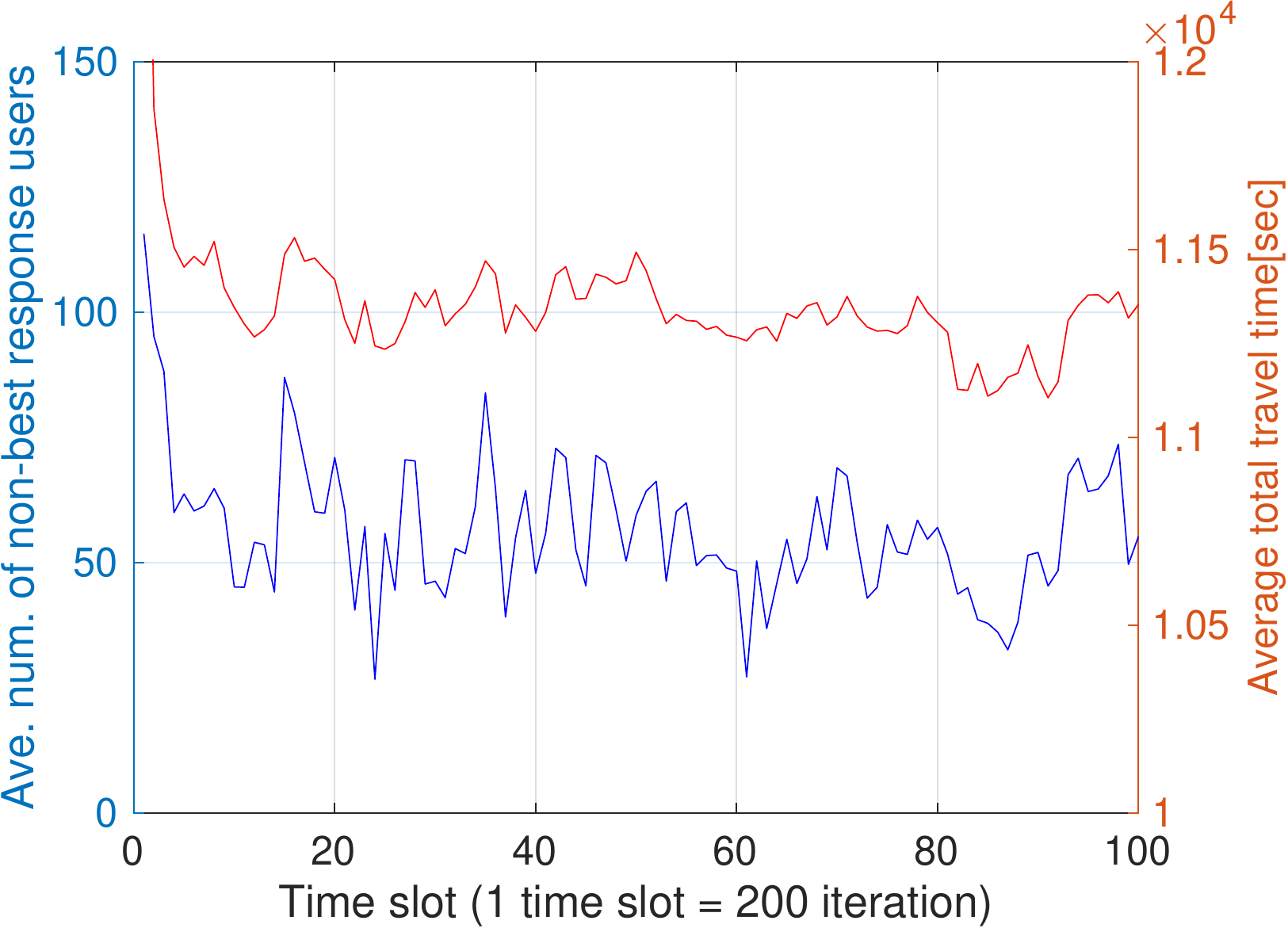}
		\subcaption{With the logarithmic perturbation parameter}\label{Fig:DSO_LogLACon}
	\end{minipage}
	\begin{minipage}[t]{0.49\textwidth}
		\centering
		\includegraphics[clip, width=0.9\columnwidth]{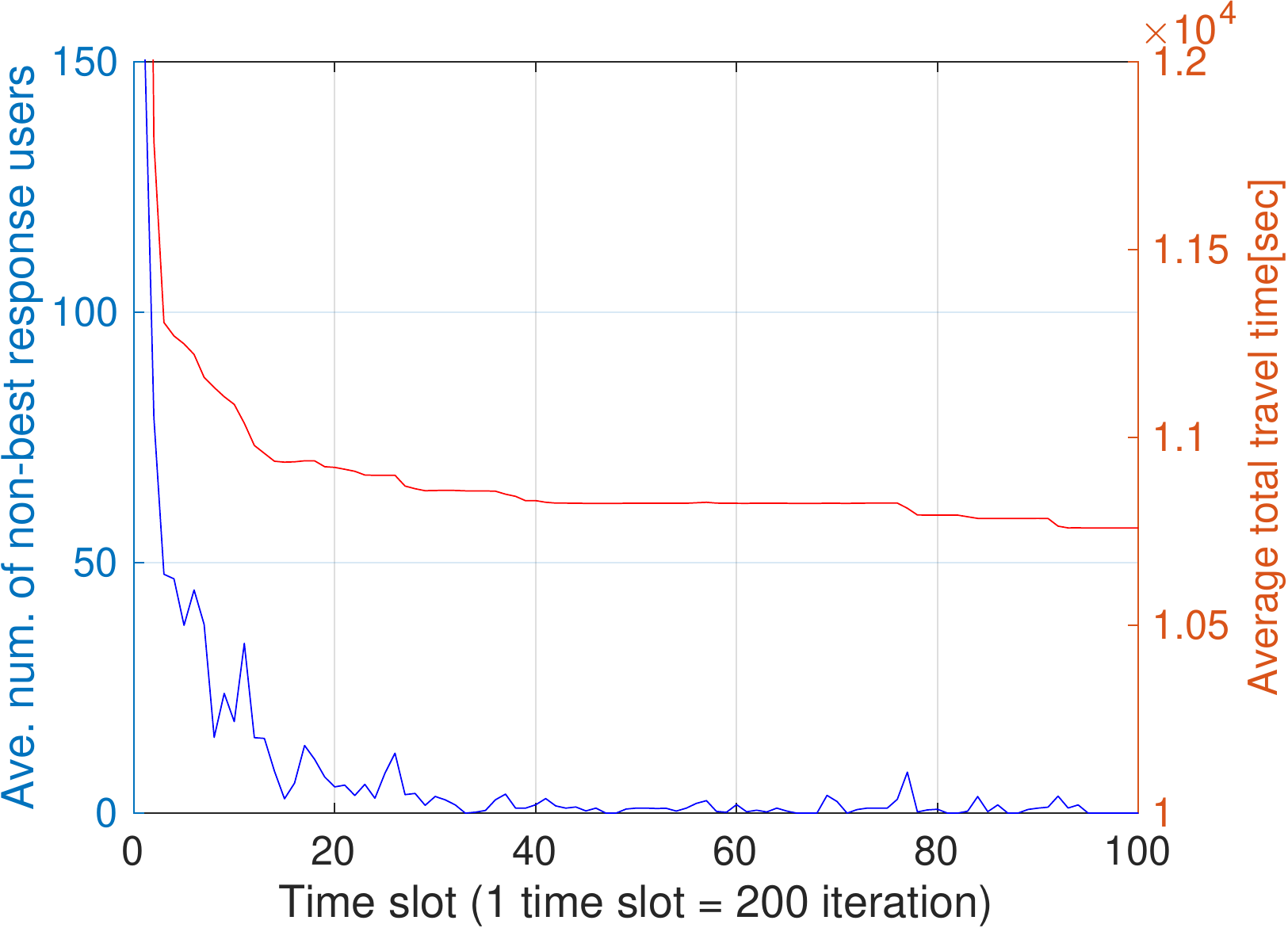}
		\subcaption{With the linear perturbation parameter}\label{Fig:DSO_LogPCon}
	\end{minipage}
	\vspace{-1mm} 
	\caption{Changes in the average number of non-best response users and \textcolor{black}{total cost} produced by the stochastic solution algorithms}
	\label{Fig:Results_S_Process}
	\vspace{-0mm}
\end{figure}

Figure~\ref{Fig:TTTDSO_S} shows the distribution of the best \textcolor{black}{total cost} obtained in each sample path of the logit dynamics. 
The result of the best response dynamics is also shown in the figure for a reference. 
We see that the stochastic algorithm with the logarithmic perturbation parameter tends to produce higher \textcolor{black}{total costs} compared to those produced by the best response dynamics. 
By contrast, the stochastic algorithm with the linear perturbation parameter tends to produce lower \textcolor{black}{total costs} although the difference in its values is relatively small. 
These results can be explained by the difference between the adjustment rules of parameter $\beta(\tau)$: the speed of the decrease in the perturbation parameter in the former case is slow, and thus the traffic states are continually perturbed during the iterations. 
\textcolor{black}{Meanwhile, traffic states under the logit response dynamics with the linear perturbation parameter are not continually perturbed owing to the fast decreasing speed of the perturbation parameter. 
However, since probabilities of deviating from inefficient optimal states become significantly small in the end, traffic states sometimes get stuck in the inefficient states in the numerical experiments with the practical number of iterations.}
The validity of this explanation is confirmed by Figure~\ref{Fig:Results_S_Process} that shows the stochastic processes (sample paths) by the above two stochastic algorithms from the same initial state.
In summary, the stochastic algorithm with the linear perturbation parameter may achieve a lower \textcolor{black}{total cost} in a practical number of iterations.

\begin{figure}[t]
	\begin{center}
	\hspace{0mm}
    \includegraphics[width=90mm,clip]{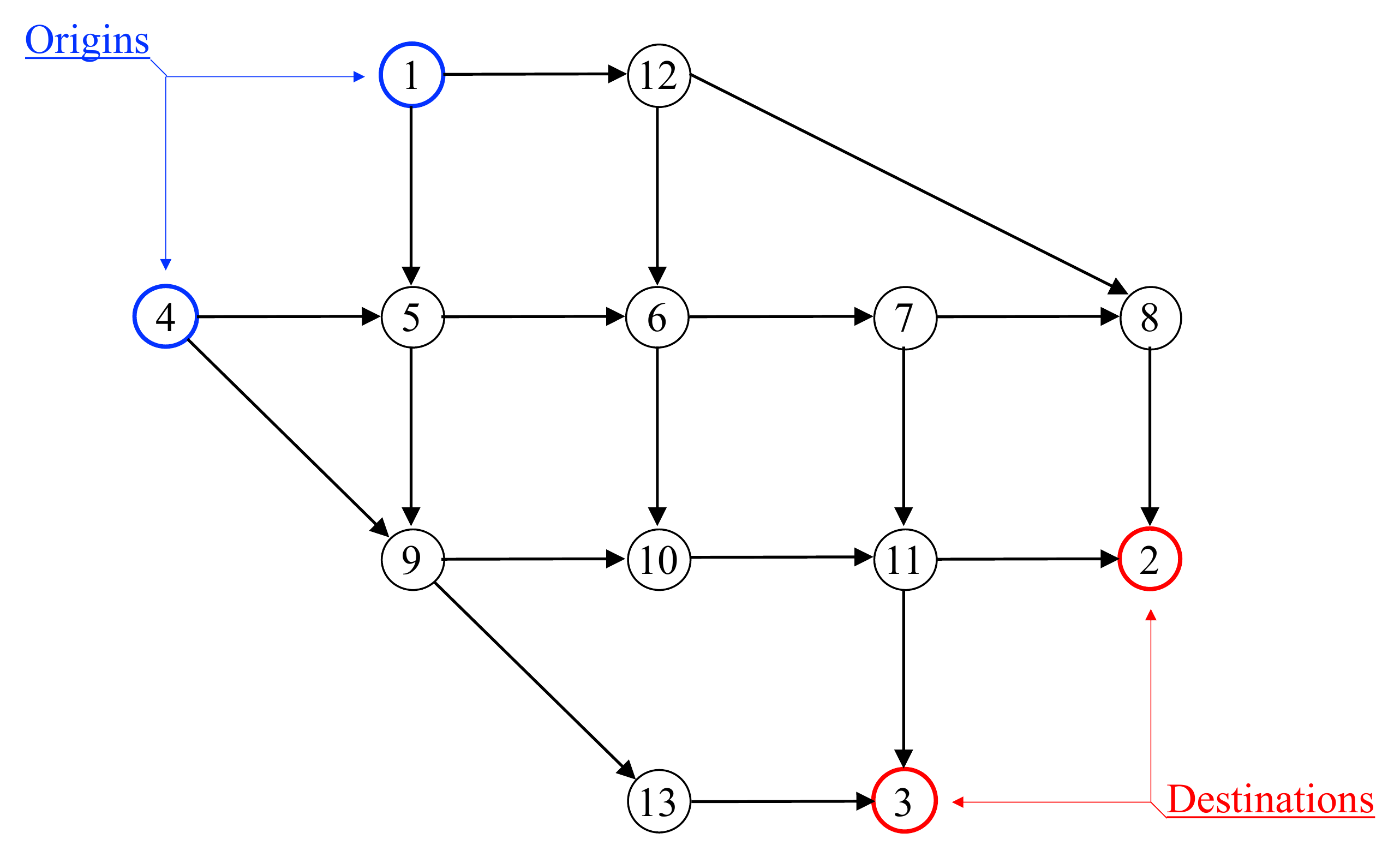}
	\end{center}
    \vspace{-4mm}
	\caption{Nguyen-Dupuis network}
    \vspace{-0mm}
    \label{Fig:NDNetwork}
\end{figure}

\begin{table}[t]
\caption{Physical conditions of each link (FFTT: Free Flow Travel Time, BC: Bottleneck Capacity, SF: Saturation Flow)}
\centering
\scalebox{0.7}{
	 \begin{tabular}{lp{1.5cm}p{1.5cm}p{1.5cm}lp{1.5cm}p{1.5cm}p{1.5cm}}\hline
	Link & FFTT [sec] & \textcolor{black}{BC} [veh/sec] & SF [veh/sec]&Link & FFTT [sec] & \textcolor{black}{BC} [veh/sec] & SF [veh/sec]\\ \hline
	$(1,5)$ 	& 		42   & 	1.25 	&	 	6 		&  	$(8,2)$ 		&  72	&	1.25		&	6\\
	$(1,12)$ 	&  	54 	& 1.25 		& 		6 		&     $(9,10)$		&  60	&	2.25		&	6\\
	$(4,5)$  	&  	54  	& 1.25 		&	 	6 		&     $(9,13)$ 		&  54	&   0.83	&	6\\
	$(4,9)$ 	&  	90	&   0.83	 &   	6  	&     $(10,11)$		&  18	&   1.67	&	6\\
	$(5,6)$ 	&  	36	&	1.42		 &		6		&     $(11,2)$		&  54	&   1.25	&	6\\
	$(5,9)$ 	&  	54	&	1.67		 & 	6 		&     $(11,3)$		&  42	&  1.25		&  6\\
	$(6,7)$ 	&  	24	&	1.25		 &		6		&     $(12,6)$		&  30	&	0.67		&	6\\
	$(6,10)$ 	&  	78	&	0.83		&		6		&  	$(12,8)$  		&  84	&		1.25	&	6\\ 
	$(7,8)$ 	&  	48	&	0.83		&		6		&		$(13,3)$		&  66	&	0.83		&	6\\ 
	$(7,11)$ 	&  	66	&	0.92		&		6		&   						&   		& 				&   \\    \hline
	 \end{tabular}
	  \label{Tab:Linkset}
	 }
\end{table}

\subsubsection{\textcolor{black}{Nguyen-Dupuis network}}

\begin{figure}[t]
	\begin{center}
	\hspace{0mm}
    \includegraphics[width=70mm,clip]{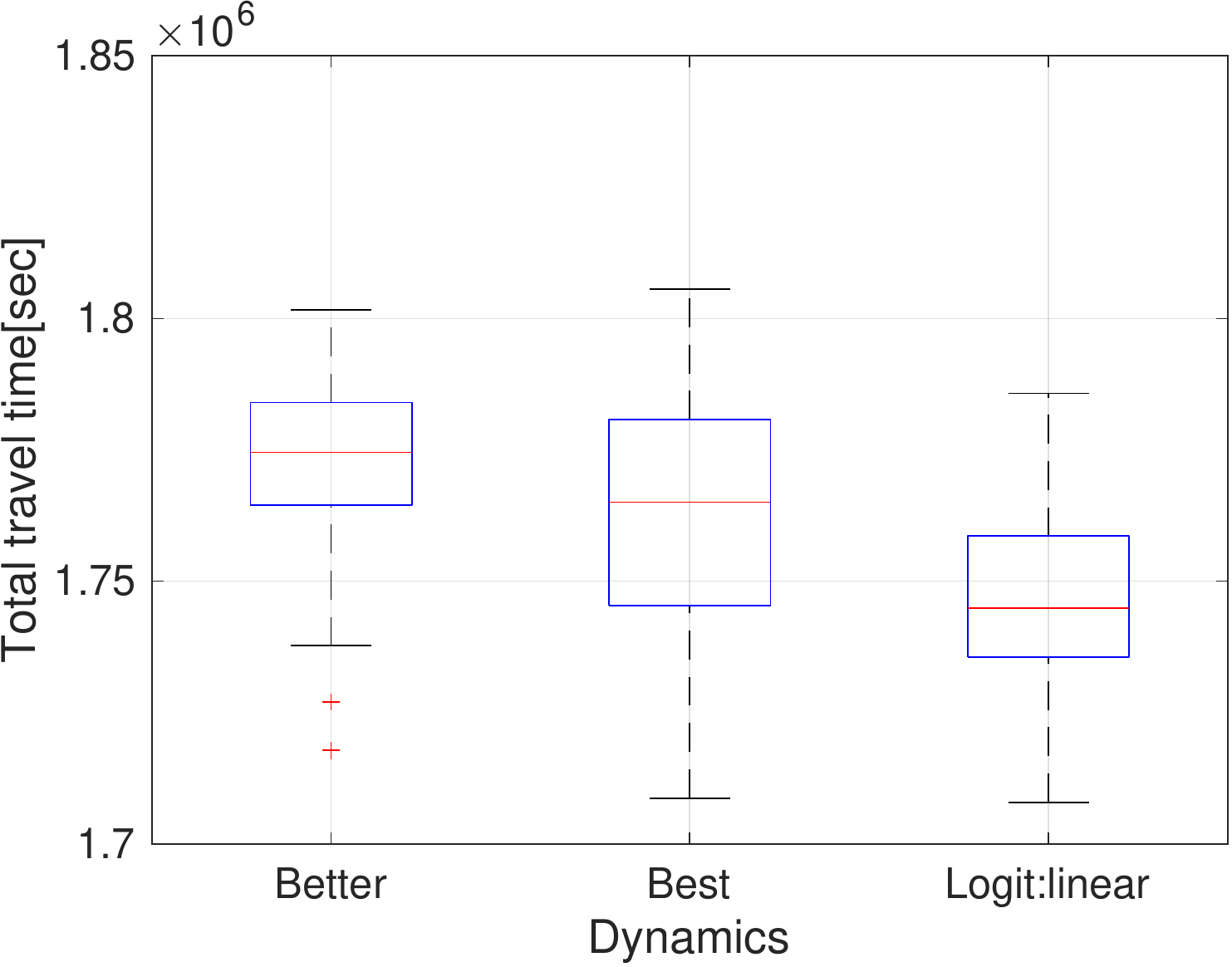}
	\end{center}
    \vspace{-4mm}
	\caption{\textcolor{black}{Total costs} of the best optimal states in sample paths in the Nguyen-Dupuis network}
    \vspace{-0mm}
    \label{Fig:Results_ND}
\end{figure}


We conduct a similar experiment in the \textcolor{black}{Nguyen-Dupuis network with multiple OD pairs (Figure~\ref{Fig:NDNetwork}). }
The physical conditions of each link are summarised in Table \ref{Tab:Linkset}.
The number of users departing from each origin-destination pair is 1,000 (i.e. 4,000 in total).
The users depart from each origin with a fixed time headway, \SI{0.5}{sec/veh}.
We generate 50 samples for each dynamics and the number iteration is set as $200,000$; we employ the logit response dynamics with a time-varying parameter $\beta(\tau) = (\tau+1)/100,000$ as the stochastic algorithm.


The results in Figure~\ref{Fig:Results_ND} show that the solution algorithms have the same qualitative properties in the previous section: the achieved \textcolor{black}{total costs} are good in the following order: the logit, best and better response dynamics.
\textcolor{black}{These results are consistent with the properties in Section~\ref{Sec:PropDSOgame}.}
Future work should investigate the characteristics of flow patterns in different optimal states.

\section{Marginal cost pricing and stochastic evolutionary implementation}\label{Sec:DUEFCP}
A direct and the most important application of DSO assignment is to provide insights into the marginal cost or Pigouvian pricing for achieving optimal states (as Nash equilibrium). 
In this context, the analyses of the evolutionary dynamics so far can be interpreted as those of the \textit{evolutionary implementation scheme of the marginal cost pricing} \citep{Sandholm2002,Sandholm2005,Sandholm2007}, in which toll level is adjusted according to the realised traffic state on a day-to-day basis.
In this section, we examine whether such an implementation scheme is essentially important for achieving optimal states. 
Specifically, we show some important features of the evolutionary implementation scheme through a comparison with another typical scheme, `fixed implementation scheme', in which toll level is set optimally in advance according to a (known or target) optimal state.

\begin{figure}[t]
	\begin{center}
	\hspace{0mm}
    \includegraphics[width=90mm,clip]{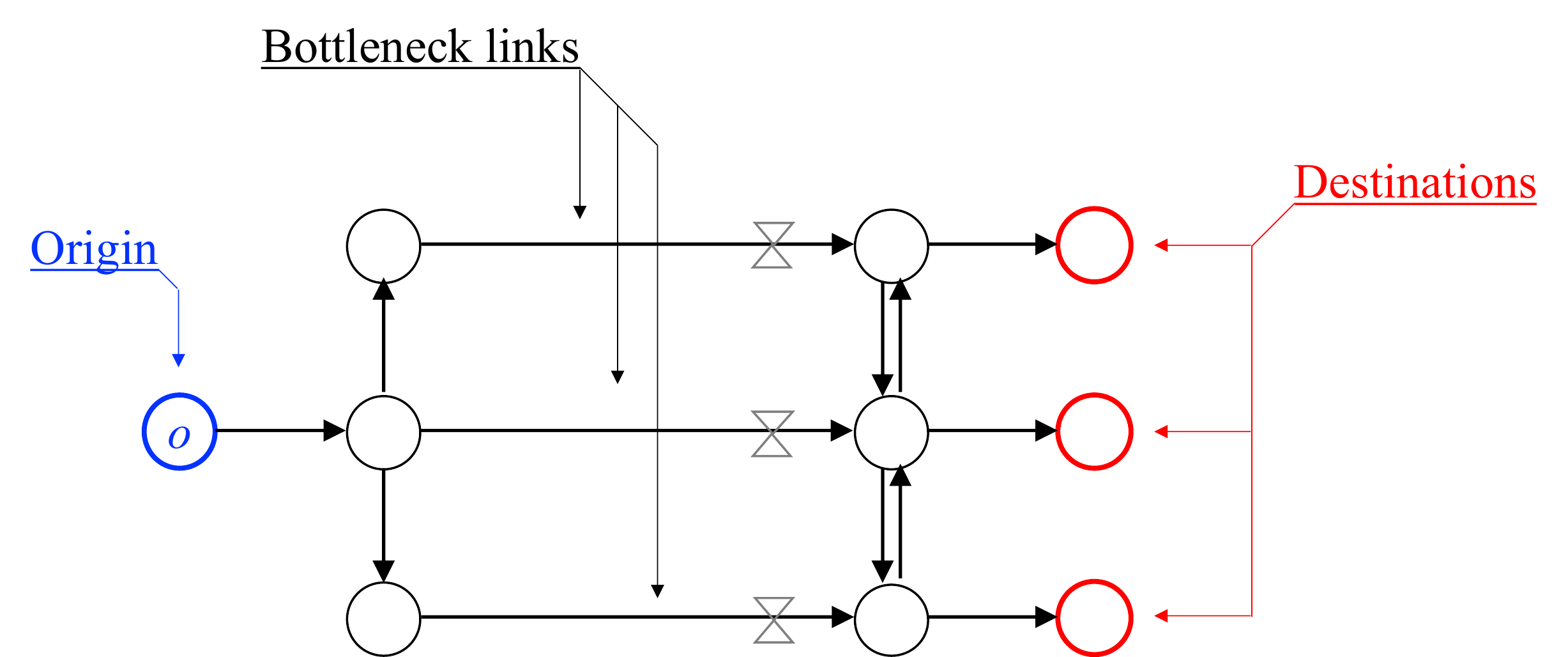}
	\end{center}
    \vspace{-4mm}
	\caption{Example of a single-bottleneck-per-route network with a single origin}
    \vspace{-0mm}
    \label{Fig:OBPRnetwork}
\end{figure}

\textcolor{black}{
This section aims to clarify differences between the two schemes from rigorous theoretical analysis.
To this end, we mainly focus on a simple network in which there exists a single origin and each route contains only one capacitated link, i.e. a single bottleneck per route network, as shown in Figure~\ref{Fig:OBPRnetwork}. 
We refer to this network as an `SBPR-1 network', in short.
SBPR-1 networks have some limitations regarding the topology.
However, these networks can be considered as sub-networks included in general networks.
Convergence and stability properties in such sub-networks will strongly affect those in whole networks as equilibrium of traffic flow of each OD pair is necessary for the equilibrium of whole traffic flow.
This implies that qualitative characteristics of the convergence and stability properties in general networks would inherit those in SBPR-1 networks.}

In Section \ref{sec:DUE-FCP} and \ref{sec:WAG}, we first show some theoretical properties of a new game under fixed tolls. 
In Section \ref{sec:Comparison_theory}, we compare the two implementation schemes theoretically.
Then, we conduct numerical experiments to validate the theoretical results, and examine whether the theoretical results are still applicable beyond an SBPR-1 network in Section \ref{sec:Comparison_Numerical}.

\subsection{DUE-FCP game and its properties}\label{sec:DUE-FCP}

We formulate a strategic game where (optimal) fixed tolls are imposed on users; we refer to this game as a `dynamic user equilibrium game with a fixed congestion pricing (DUE-FCP game)'.
\textcolor{black}{In the DUE-FCP game, we assume that the utility of each user is equal to the negative value of the sum of the route travel time of the user and fixed toll.
Mathematically, the utility of user $i$ who selects the route $r_{i}\in\mathcal{R}_{i}$ for a given route profile $\mathbf{r}_{-i}$ is represented by,}
\begin{align}
&U^{F}_{i}(r_{i}, \mathbf{r}_{-i}) = - C_{i}(r_{i}, \mathbf{r}_{-i}) - T_{i}(r_{i}), \label{Eq:Utility_Fixed}
\end{align}
where $T_{i}(r)$ is the toll imposed on user $i$ whose route is $r$.
A Nash equilibrium state $\mathbf{r}^{*}$ is then defined as a state where users cannot improve their utility by unilaterally changing routes, as follows:
\begin{align}
U^{F}_{i}(r^{*}_{i}, \mathbf{r}^{*}_{-i}) = \max_{r\in\mathcal{R}_{i}}U^{F}_{i}(r, \mathbf{r}^{*}_{-i}),\quad \forall i\in\mathcal{P}.\label{Eq:Nash_Fixed}
\end{align}

In this game, with a fixed set of tolls, one can make a certain traffic state a Nash equilibrium state, as follows:

\begin{prop}\label{Prop:FixedTollAdjust}
Consider a DUE-FCP game and an arbitrary route profile $\mathbf{r}^{*}\in \mathcal{R}$.
Then, there always exists at least one set of tolls $\mathbf{T}^{*}$ that makes $\mathbf{r}^{*}$ into strict Nash equilibrium.
That is, the following condition is satisfied under $\mathbf{T}^{*}$:
\begin{align}
U^{F}_{i}(r^{*}_{i}, \mathbf{r}^{*}_{-i}) > U^{F}_{i}(r, \mathbf{r}^{*}_{-i}),\quad \forall r\in\mathcal{R}_{i}\setminus \{  r^{*}_{i} \},\forall i\in\mathcal{P}.\label{Eq:TollAdjustment}
\end{align}
\end{prop}
\begin{prf}
It is obvious that changes in tolls of routes for a user do not change the utility of the other users.
Thus, without affecting the other users' utilities, we can make $r^{*}_{i}$ the unique best response route for user $i$ in $\mathbf{r}^{*}$ by setting tolls of each route to $T^{*}_{i}(r)$ as follows: 
\begin{align}
&- C_{i}(r^{*}_{i}, \mathbf{r}^{*}_{-i}) - T^{*}_{i}(r^{*}_{i}) > 
- C_{i}(r_{i}, \mathbf{r}^{*}_{-i}) - T^{*}_{i}(r_{i})   \quad \forall  r_{i}\in\mathcal{R}_{i}\setminus\{r_{i}^{*} \}
\end{align}
By applying the same procedure to all the users, we ensure that $\mathbf{r}^{*}$ is strict Nash equilibrium in the DUE-FCP game. \qed
\end{prf}
\noindent If the tolls for each user are set as the negative externalities that he/she creates in an optimal state $\mathbf{r}^{*}$ (a Nash equilibrium state in a DSO game), i.e. $T_{i}(r_{i}) = E_{i}(r_{i},\mathbf{r}^{*}_{-i})$, a Nash equilibrium state in the DUE-FCP game corresponds to the optimal state. 
Thus, a DUE-FCP game with a set of optimal tolls is useful for analysing the properties of the fixed implementation scheme. 
\textcolor{black}{From now on, we focus on an SBPR-1 network as mentioned above.}
This network has the following \textit{ordering property} that characterises the direction of the externalities: 
\begin{prop}\label{Prop:OrderingProperty}
Consider a DUE-FCP game in an SBPR-1 network.
The utility of any route of a user is independent of the route choices of the users who depart from the origin later than the considered user.
\end{prop}
\begin{prf}
The utility of a user's route consists of the route travel time and corresponding fixed toll.
As the travel time in every uncapacitated link does not change, it is sufficient for us to prove that the travel time of a user in each capacitated link is independent of the route choices of the users who depart later than the user.

As no delay occurs in the uncapacitated links, a user arrives at any node that is situated upstream of the capacitated links earlier than users who depart later than the user.
\textcolor{black}{Note that this holds even if perturbations exist because the actual travel time from the origin to the capacitated link selected by the user does not change.}
By combining this fact with the FIFO principle and the causality of the dynamic loading model, it is guaranteed that the travel time of a capacitated link, which includes delays owing to traffic congestion, is independent of the route choices of the users who depart later than the user.
In other words, it is dependent only on the route choices of users who depart earlier than the user. \qed
\end{prf}

\noindent With this ordering property, we can prove the uniqueness of equilibrium in the DUE-FCP game:
\begin{coro}\label{Coro:FixedTollAdjust}
Consider a DUE-FCP game in an SBPR-1 network and an arbitrary route profile $\mathbf{r}^{*}\in \mathcal{R}$.
Consider also a set of tolls $\mathbf{T}^{*}$ that makes $\mathbf{r}^{*}$ into strict Nash equilibrium.
Then, $\mathbf{r}^{*}$ becomes the unique Nash equilibrium in the DUE-FCP game.
\end{coro}
\begin{prf}
See \ref{App:FixedTollAdjust}. \qed
\end{prf}

\noindent The ordering property is useful for analytically investigating the properties of equilibrium, as utilised in previous studies~\citep[e.g.][]{Kuwahara1993,Akamatsu2015,Wada2019}.
Especially, with this property and the theory of weakly acyclic games, \cite{Satsukawa2019} established the convergence and stability of dynamic user equilibrium in a unidirectional network. 
In a similar manner, we prove the convergence and stochastic stability of evolutionary dynamics in a DUE-FCP game in an SBPR-1 network.

\subsection{Weakly acyclicity, convergence and stability of DUE-FCP games in an SBPR-1 network}\label{sec:WAG}

To define a weakly acyclic game, we utilise the notion of a better response path.
A better response path is a sequence of strategy profiles $\mathbf{r}^{1}, \mathbf{r}^{2},\ldots, \mathbf{r}^{L}$ such that for each successive pair $\mathbf{r}^{\tau}, \mathbf{r}^{\tau+1}$, there is exactly one user $i^{\tau}$ that changes his/her route and that user improves the utility, as follows:
\begin{align}
\begin{cases}
r^{\tau}_{i^{\tau}}  \neq r^{\tau+1}_{i^{\tau}},  \quad \text{s.t.} \quad &r^{\tau+1}_{i^{\tau}}\in D_{i^{\tau}}(\mathbf{r}^{\tau}) ,\\
r^{\tau}_{j}  \neq r^{\tau+1}_{j} ,\quad    &j\in\mathcal{P}\setminus \{  i^{\tau} \}.
\end{cases}
\end{align}

\noindent The class of weakly acyclic games, introduced by \cite{Young1993} and \cite{Marden2012b}, is then defined as follows:
\begin{defi}
A game is a \textit{weakly acyclic game} if from any strategy profile $\mathbf{r}\in\mathcal{R}$, there exists a better response path starting at $\mathbf{r}$ and ending at a Nash equilibrium state.
\end{defi}
We then obtain the following theorem and propositions regarding theoretical properties of a DUE-FCP game in an SBPR-1 network in the same way as in \cite{Satsukawa2019}:

\begin{theo}\label{Theo:DUE-FCPisWAG}
A DUE-FCP game in an SBPR-1 network is a weakly acyclic game.
\end{theo}
\begin{prf}
The outline of the proof is described as follows.
Thanks to \textbf{Proposition~\ref{Prop:OrderingProperty}}, it is guaranteed that if a user selects a best response route in a particular route profile, the route remains a best response route regardless of any route changes made by the users departing later.
Thus, by changing the routes of the users to their best response routes sequentially in the order of their departure times, all users come to select their ex-post best response routes, that is, we obtain a Nash equilibrium state (see \ref{App:DUE-FCPisWAG} for details). \qed
\end{prf}


\begin{prop}
Consider a DUE-FCP game in an SBPR-1 network, an arbitrary route profile $\mathbf{r}^{*}\in \mathcal{R}$, and a set of tolls $\mathbf{T}^{*}$ that makes $\mathbf{r}^{*}$ into a unique strict Nash equilibrium.
Then, a route profile generated by the better/best response dynamics from an arbitrary initial route profile converges almost surely to $\mathbf{r}^{*}$, that is, $\mathbf{r}^{*}$ is global asymptotically stable under the best response dynamics.
\end{prop}
\begin{prf}
\cite{Young2004} proved the convergence of the better response dynamics to Nash equilibrium in weakly acyclic games.
Since $\mathbf{r}^{*}$ is the unique Nash equilibrium, the global asymptotic stability is promptly derived. 
Next, in the same way as the proof of \textbf{Theorem~\ref{Theo:DUE-FCPisWAG}}, we can prove that there exists a sequence of best responses (i.e. best response path) from any initial route profile ending at the strict Nash equilibrium state.
This game is called a weakly acyclic game under best replies, and \cite{Young2004} proved the convergence of the best response dynamics in the game: thus, the proposition holds.
\qed
\end{prf}

\begin{prop}\label{Prop:SSS-DUEFCP}
Consider a DUE-FCP game in an SBPR-1 network, an arbitrary route profile $\mathbf{r}^{*}\in \mathcal{R}$, and a set of tolls $\mathbf{T}^{*}$ that makes $\mathbf{r}^{*}$ into a unique strict Nash equilibrium.
Then, $\mathbf{r}^{*}$ is a stochastically stable under the logit response dynamics.
\end{prop}
\begin{prf}
This is promptly derived by the Theorem 4 of \cite{Young1993}.\qed
\end{prf}

\subsection{Theoretical comparison: Evolutionary vs. fixed implementation schemes}\label{sec:Comparison_theory}
\textcolor{black}{We are now ready to discuss the structural differences between the evolutionary and fixed implementation schemes of the marginal cost pricing in SBPR-1 networks.}
The convergence and stability results for the DSO and DUE-FCP games so far show that an optimal traffic state will be achieved by natural evolutionary dynamics under both implementation schemes.
However, there are two major differences in the mechanism for achieving optimal states between the two schemes.

The first difference stems from the differences in externalities that are created under each implementation scheme, which will affect convergence processes. 
Under the fixed implementation scheme (or in a DUE-FCP game), owing to the ordering property shown in \textbf{Proposition~\ref{Prop:OrderingProperty}}, a route change of a user affects the utility of users who depart later from the same origin, but not \textit{vice versa}.
By this temporal asymmetry of the interaction, users having earlier departure times become likely to choose their \textit{ex-post} best response routes earlier during the convergence process \citep{Satsukawa2019}. 
By contrast, under the evolutionary implementation scheme (or in a DSO game), as externalities are perfectly internalised, the better/best response of \textit{any user} leads the traffic state toward an equilibrium state, regardless of the order of his/her departure time.    
This implies that the total travel time decreases smoother to an efficient state under the evolutionary implementation scheme.



The second difference is related to the obvious difference in the price-setting, i.e., the toll level is adjusted according to the realised state on a day-to-day basis under the evolutionary implementation scheme, while the level is set in advance according to a \textit{known} target state under the fixed implementation scheme. 
This difference affects which state is stabilised.  
Under the evolutionary implementation scheme, the most efficient equilibrium state is stochastically stable, as shown in \textbf{Theorem~\ref{Theo:SSS-DSO}}.
In other words, that equilibrium state will be achieved thanks to \textit{perturbations} even though there are multiple equilibria (i.e. locally optimal states). 
On the other hand, under the fixed implementation scheme, the target state is a unique and stochastically stable equilibrium state, as shown in \textbf{Proposition~\ref{Prop:SSS-DUEFCP}}\footnote{\textcolor{black}{For static traffic assignment problems dealing with continuum users, \cite{Sandholm2002,Sandholm2005,Sandholm2007} analysed the problems with marginal cost pricing schemes by a similar approach.
The author showed that static traffic assignment problems under evolutionary and fixed implementation schemes are potential games, and obtained similar conclusions about the convergence and stability properties.
These papers also provide an excellent discussion on information issues in the implementation of congestion pricing schemes.}}.  
It means that (i) a locally (or inefficient) optimal state is stabilised unless information regarding a globally optimal state is known in advance, and (ii) perturbations do not play an essential role, unlike the evolutionary implementation scheme. 
These differences imply that when such an inefficient optimal state is a target state under the fixed implementation scheme, a traffic state could not reach a more efficient traffic state, such as globally optimal state minimising the total cost.

\subsection{Numerical experiments: Evolutionary vs. fixed implementation schemes}\label{sec:Comparison_Numerical}
From the above theoretical discussions, we can expect that under the evolutionary implementation scheme, the traffic state decreases smoother to an efficient state and a traffic state reaches a more efficient traffic state.
In this section, we conduct numerical experiments for two purposes: (i) confirming the validity of this expectation and (ii) testing whether the theoretical results (especially, for a fixed implementation scheme) are still applicable beyond an SBPR-1 network. 

\subsubsection{Simple network}
We consider the simple network used in the previous section (Figure~\ref{Fig:Network}). 
This network is an SBPR-1 network. 
To simulate the evolutionary implementation scheme, we iterate the logit response dynamics with a fixed perturbation parameter in the DSO game. 
Regarding the fixed implementation scheme, we must determine a target state and the (optimal) fixed tolls for implementing that state first.
In the experiment, we set an equilibrium route profile ${\bf r}^{*}$ obtained under the evolutionary scheme as the target state; the external costs at the state are set as fixed tolls (i.e., $T_{i}(r) = E_{i}(r, {\bf r}^{*}_{-i})$). 
We then iterate the logit response dynamics in the DUE-FCP game from the same initial state as with the DSO game.

\begin{figure}[t]
	\begin{center}
	\hspace{0mm}
    \includegraphics[width=80mm,clip]{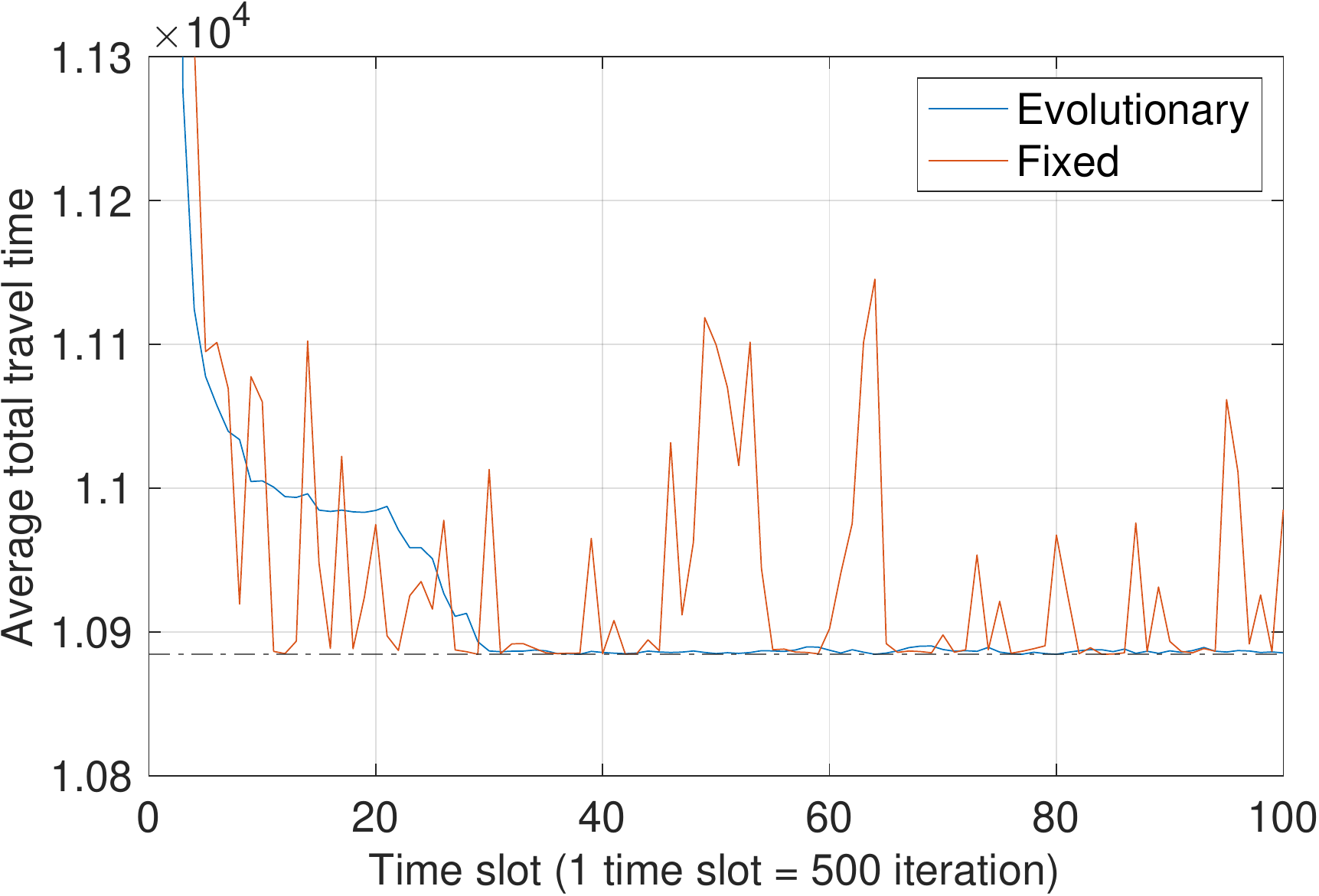}
	\end{center}
    \vspace{-4mm}
	\caption{Behaviour of \textcolor{black}{total costs} under evolutionary and fixed implementation schemes in the simple network}
    \vspace{-3mm}
    \label{Fig:DSODUELogitTTT}
\end{figure}

First, let us look at a sample path of the stochastic process under each implementation scheme.
We here set \textcolor{black}{$\beta = 1.0$}, the number of iterations as 50,000 times, and the initial route profile of each user as the shortest-distance route; the target state of the fixed implementation scheme is an efficient equilibrium state that has minimum \textcolor{black}{total cost} among equilibrium states that occurred under the evolutionary implementation scheme. 
Figure~\ref{Fig:DSODUELogitTTT} shows that the process of the average \textcolor{black}{total cost} which is calculated for every time slot (1 time slot = 500 iterations).
The black dotted line indicates the \textcolor{black}{total cost} at the efficient equilibrium.
\textcolor{black}{From this figure, we first confirm that traffic states sometimes deviate from efficient states by perturbations under the two implementation schemes.
We also see that:}
\begin{itemize}
    \item Under the fixed implementation scheme, although the process approaches the efficient equilibrium (or target) state somewhat quickly, it does not stay at the state, and frequently deteriorates.
    \item Under the evolutionary implementation scheme, the process does not significantly deteriorate during the iterations; the \textcolor{black}{total cost} decreases almost monotonically.
\end{itemize}
These results suggest that a traffic state under the evolutionary implementation scheme is more robust to perturbations than that under the fixed one in the sense that the \textcolor{black}{total costs} do not largely oscillate.


The above difference can be attributed to the differences in externalities between the two implementation schemes discussed in the previous subsection.
Under the evolutionary implementation scheme, as the utility of each user is aligned with the \textcolor{black}{total cost}, if one user makes a mistake and the resulting \textcolor{black}{total cost} increases, any user who can improve the utility rationally chooses a route so that the \textcolor{black}{total cost} decreases.
Thus, \textcolor{black}{total costs} will not significantly deteriorate by perturbations.
On the other hand, under the fixed implementation scheme, an improvement in the utility of a user does not necessarily decrease the \textcolor{black}{total cost}.
Moreover, once a traffic state is perturbed, the ordering property requires the sequential route changes of users to the best response routes in the order of their departure times to return the traffic state.
This would prevent the traffic states from recovering to an efficient state smoothly.
Thus, we can infer from these facts that the \textcolor{black}{total cost} under the evolutionary implementation scheme would become more robust to perturbations.

\begin{figure}
	\begin{minipage}[t]{0.49\textwidth}
		\centering
		\includegraphics[clip, width=0.9\columnwidth]{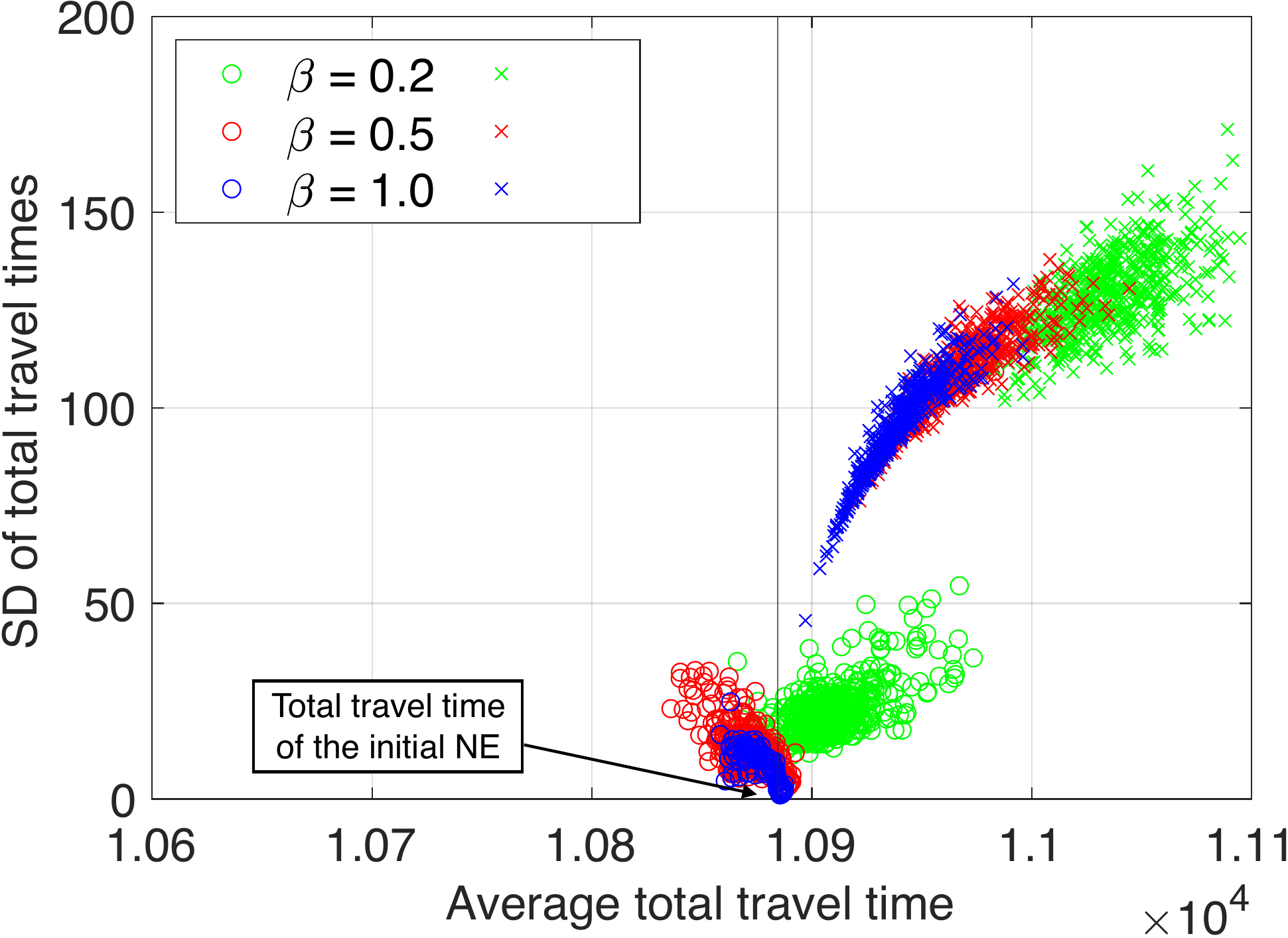}
		\subcaption{Initial state: the efficient equilibrium state in Figure~\ref{Fig:DSODUELogitTTT}}\label{Fig:TemMinNEMap}
	\end{minipage}
	\begin{minipage}[t]{0.49\textwidth}
		\centering
		\includegraphics[clip, width=0.9\columnwidth]{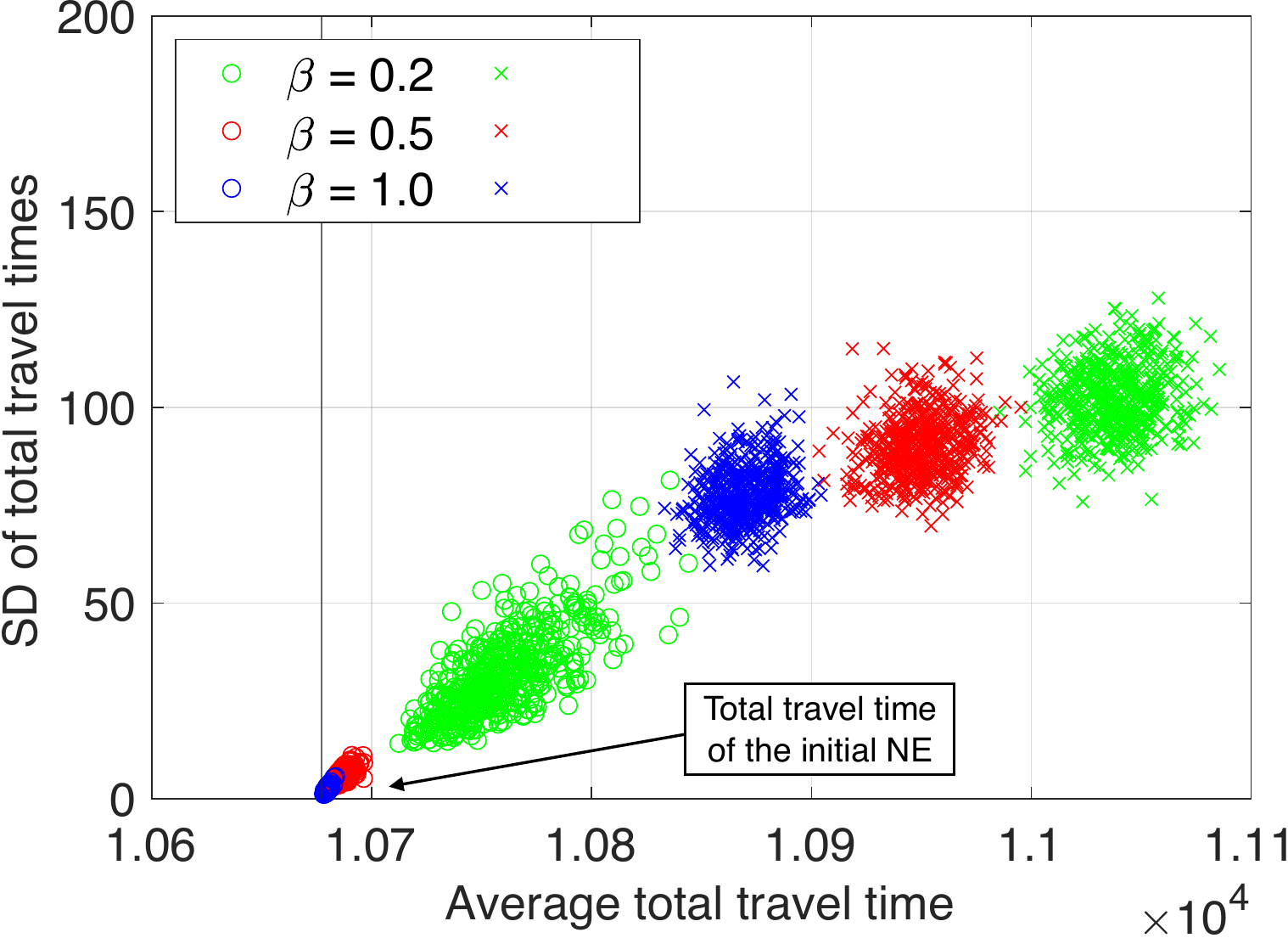}
		\subcaption{Initial state: the tentative globally optimal state.}\label{Fig:MinNEMap}
	\end{minipage}
	\vspace{-1mm} 
	\caption{Averages and standard deviations of the \textcolor{black}{total costs} around Nash equilibrium}
	\label{Fig:NEMap_Simple}
	\vspace{-0mm}
\end{figure}

To confirm the validity of the observation about the robustness, we generate 500 sample paths from an equilibrium state under each implementation scheme for every experimental setting; the number of iteration for each sample path is set as 20,000.  
\textcolor{black}{We select three perturbation parameters considering the frequencies of mistakes (i.e. changing to routes with worse utilities): $\beta = 1.0, 0.5$ and $0.2$}; and two initial states: the efficient equilibrium in Figure~\ref{Fig:DSODUELogitTTT} and a tentative globally optimal state obtained in the experiment in Section~\ref{Sec:Solutions}. 
The target state of the fixed implementation scheme is set as the initial state. 


Figure~\ref{Fig:NEMap_Simple} shows the results. 
Each mark represents the average and standard deviation of the \textcolor{black}{total costs} derived in one sample path (circle: the evolutionary implementation scheme; crosses: the fixed implementation scheme); different colours of marks indicate different perturbation parameters; each black line represents the \textcolor{black}{total cost} at each initial equilibrium.
\textcolor{black}{In Figure~\ref{Fig:TemMinNEMap}, the frequencies of mistakes are about 0.5, 1.0 and 2.0\% when $\beta$ is 1.0, 0.5 and 0.2, respectively.
On the other hand, in Figure~\ref{Fig:MinNEMap}, the frequencies are about 0.5, 1.0 and 2.0\% under the evolutionary implementation scheme, and 1.0, 1.5 and 2.5\% under the fixed implementation scheme.}

Figure~\ref{Fig:TemMinNEMap} shows that the standard deviations under the fixed implementation scheme tend to be larger than those under the evolutionary one; this tendency is particularly evident when the perturbation level becomes high (from \textcolor{black}{$\beta = 1.0$ to $\beta = 0.2$}). 
The tendencies are also observable in Figure~\ref{Fig:MinNEMap}, which shows the results of the process starting from the tentative globally optimal state. 
These results are consistent with the above observations that the \textcolor{black}{total cost} under the evolutionary implementation scheme is more robust to perturbations than the fixed implementation scheme.  

\textcolor{black}{Figures~\ref{Fig:TemMinNEMap} and~\ref{Fig:MinNEMap} also show that the average total costs under the fixed implementation scheme are systematically higher than the initial one.
This is because the (initial) locally optimal state becomes strongly stabilised by the fixed tolls and a traffic state is attracted to the locally optimal state.
Owing to this, the traffic state tends to stay among those with higher total travel times near the locally optimal state.
Thus, the total travel time becomes worse in average than the initial one.
Meanwhile, some average total costs under the evolutionary implementation scheme in Figures~\ref{Fig:TemMinNEMap} are less than the total cost of the initial one.
This is because the globally optimal state is stochastically stable, and a traffic state is attracted to the most efficient state.
This observation is consistent with the theoretical discussion in the previous subsection. 
}

\begin{figure}
	\begin{minipage}[t]{0.49\textwidth}
		\centering
		\includegraphics[clip, width=0.9\columnwidth]{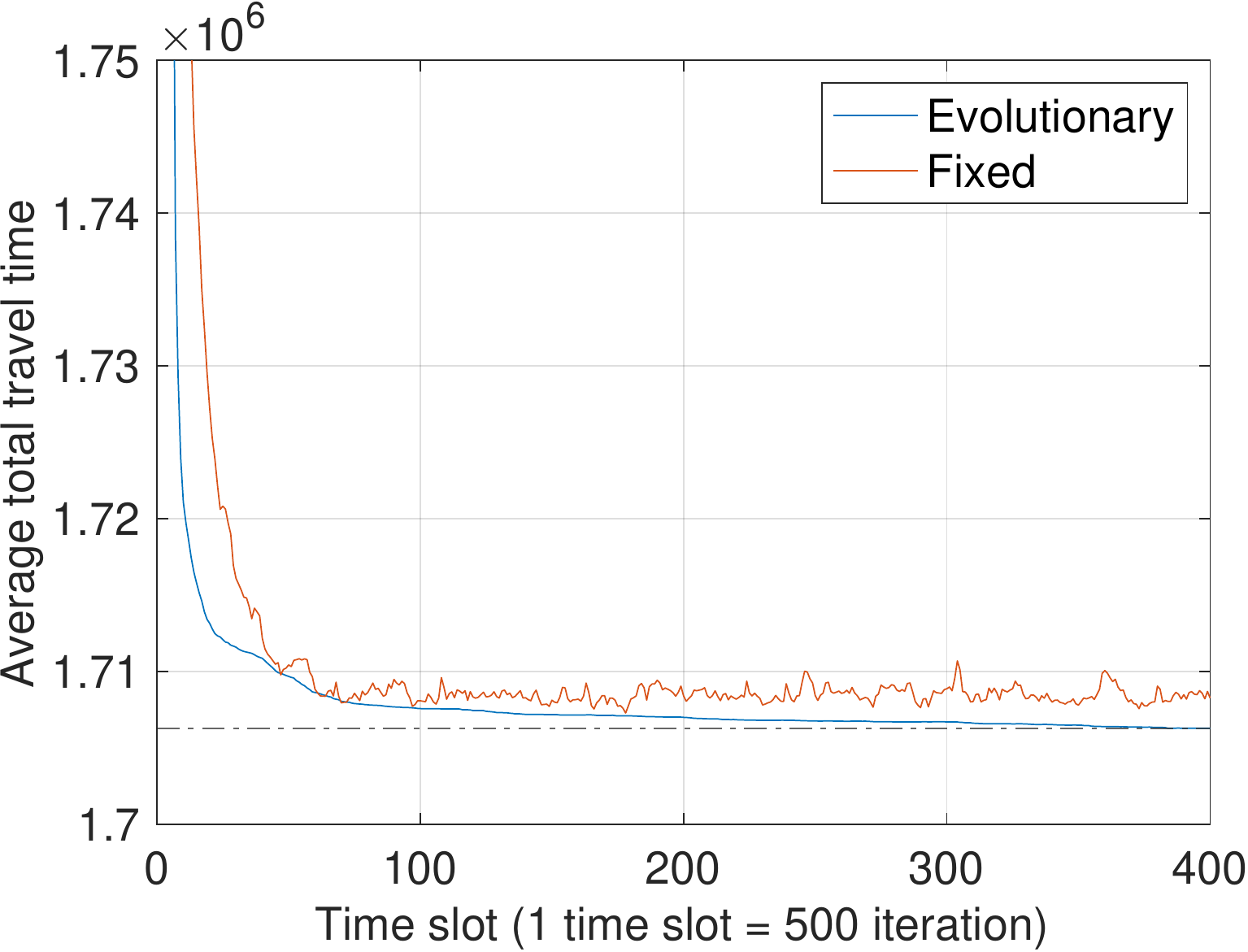}
		\subcaption{Behaviour of \textcolor{black}{total costs} from the shortest-distance route profile}\label{Fig:DSODUEAveTTT-ND}
	\end{minipage}
	\begin{minipage}[t]{0.49\textwidth}
		\centering
		\includegraphics[clip, width=0.9\columnwidth]{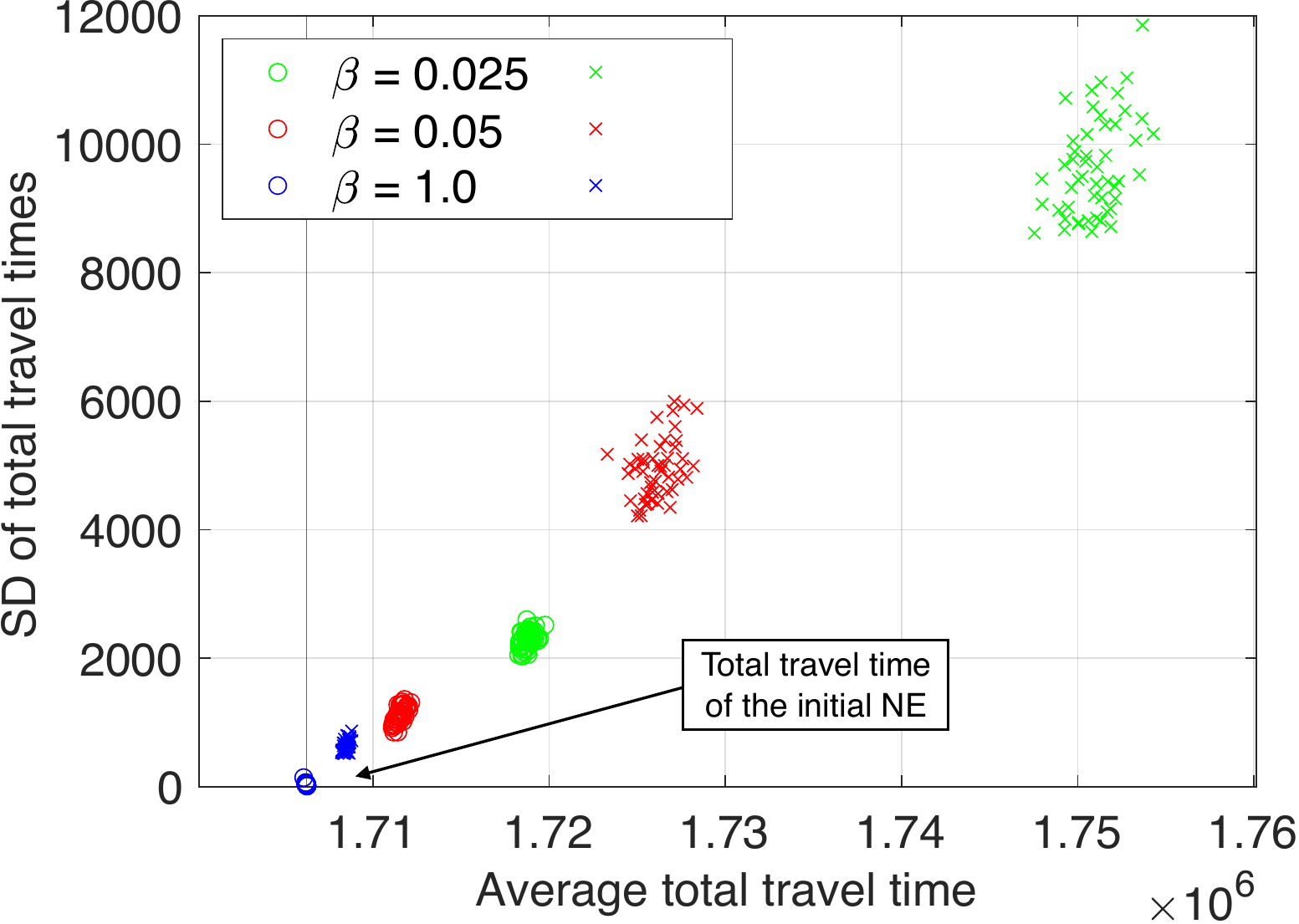}
		\subcaption{Averages and standard deviations of \textcolor{black}{total costs} around Nash equilibrium}\label{Fig:DSODUEMap-ND}
	\end{minipage}
	\vspace{-1mm} 
	\caption{Results of the behaviour of average \textcolor{black}{total costs} under the logit response dynamics in the Nguyen-Dupuis network}
	\label{Fig:DSODUEinNDnet}
	\vspace{-0mm}
\end{figure}

\subsubsection{Nguyen-Dupuis network}

Finally, we consider the \textcolor{black}{Nguyen-Dupuis network with tandem bottlenecks and multiple OD pairs (i.e. non-SBPR-1 network)}.
Figure~\ref{Fig:DSODUEAveTTT-ND} shows the sample paths of the average \textcolor{black}{total costs} (here, 1 time slot = 500 iterations) under the evolutionary and fixed implementation schemes from the shortest-distance route profile.
We set \textcolor{black}{$\beta = 1.0$, the number of iterations as 200,000 times, and iterate the dynamics from the shortest-distance route profile. 
The target state of the fixed implementation scheme is Nash equilibrium near the most efficient state with the minimum total cost under the evolutionary implementation scheme: we derived this Nash equilibrium by employing the better response dynamics where the initial route profile was set to the most efficient state.}
From this figure, we can see that the average \textcolor{black}{total costs} under the evolutionary implementation scheme decrease almost monotonically, whereas, those under the fixed implementation scheme frequently oscillate, as with the previous experiment.


Figure~\ref{Fig:DSODUEMap-ND} shows the behaviour of the average \textcolor{black}{total costs} around Nash equilibrium.
\textcolor{black}{We consider three perturbation parameters, $\beta = 1.0$, $0.05$ and $0.025$, and generate 50 sample paths from an equilibrium state under each scheme.
Note that the frequencies of mistakes are about $3$, $10$ and $15$\% when $\beta = 1.0$, $0.05$ and $0.025$, respectively.}
The number of iterations for each sample path is set as 50,000 times.
We set the most efficient Nash equilibrium as the initial state.
This figure confirms that the characteristics of the behaviour found in the SBPR-1 network also hold in a non-SBPR-1 network, also.
Specifically, the standard deviations of the \textcolor{black}{total costs} under the fixed implementation scheme are larger than those under the evolutionary implementation scheme, and the average \textcolor{black}{total costs} under the evolutionary implementation scheme tend to be lower than those under the fixed implementation scheme.
To explore whether this result holds more generally, a systematic set of numerical experiments should be conducted for different types of network settings.

\section{Conclusion}\label{Sec:Conclusion}
In this study, we examined the convergence and stability of DSO traffic assignment with atomic users.
We first formulated a DSO game that is a strategic game form of the DSO traffic assignment wherein atomic users select routes minimising their marginal social costs.
This game becomes a potential game whose potential function is the negative of \textcolor{black}{total cost function}; the Nash equilibrium states correspond to the optimal states of a problem minimising the \textcolor{black}{total cost}.
These enable us to rigorously analyse the theoretical properties of evolutionary dynamics.
We proved that the better/best response dynamics converges to a set of optimal states; a globally optimal state is stochastically stable under the logit response dynamics.
The numerical experiments demonstrated that the \textcolor{black}{total costs} generated by the best response dynamics tend to be lower than those by the better response dynamics, and the stochastic algorithm with the linear perturbation parameter might achieve a lower \textcolor{black}{total cost} in a practical number of iterations.

As an application of the analysis above, we further examined the characteristics of the evolutionary implementation scheme of marginal cost pricing by comparing this scheme with a fixed implementation scheme.
Specifically, we first showed the convergence and stability in a DUE-FCP game, which is a strategic game under the fixed implementation scheme, in an SBPR-1 network with the theory of weakly acyclic games.
From the differences in externalities and stable states between DSO and DUE-FCP games, we clarified that, under the evolutionary implementation scheme, \textcolor{black}{(i) the total cost decreases smoother to an efficient traffic state, and (ii) a traffic state could reach a more efficient state in principle.}
Finally, we conducted numerical experiments to validate the theoretical findings. 
The experiments also indicate the differences in the robustness to perturbations of the two implementation schemes: under the fixed implementation scheme, \textcolor{black}{total costs} frequently deteriorate by perturbations; in contrast, under the evolutionary implementation scheme, \textcolor{black}{total costs} decreases almost monotonically.

As part of future works, it is important to clarify traffic flow patterns of optimal states through theoretical analysis and systematic numerical experiments.
Especially regarding numerical experiments, although we provide algorithms that are guaranteed to converge, they may be inefficient and their computational costs would be high in large-scale networks.
It is thus required to develop efficient algorithms by utilising useful properties of DSO game, such as the existence of potential functions.
In relation to this, it is important to examine what types of individuals' learning behaviour and available information converge evolutionary dynamics to optimal states quickly.
To explore the applicability of the proposed approach to more general DSO assignment with both route and departure time choices is also important.
\textcolor{black}{Furthermore, it is interesting to pursue a realistic scheme that (approximately) implements the marginal cost pricing in this study, for example, by considering the link-based decomposition of the marginal cost;} 
it is also interesting to purse another implementation scheme to achieve an optimal state, such as controlling users via an autonomous vehicle system.
The DSO game could be directly applied to the case in which all vehicles are autonomous vehicles.
However, there will exist more complex situations, such as the existence of non-autonomous vehicles.
Introducing such conditions into a strategic game and analysing the theoretical properties are important topics for future studies.

\section*{Acknowledgements}
This work was financially supported by JSPS KAKENHI Grant numbers JP20K14843 and JP20H00265.
The comments of anonymous reviewers are gratefully acknowledged.

\appendix

\section{A dynamic loading model dealing with atomic users}\label{App:DNL}

\subsection{Link model}
\begin{figure}[t]
	\begin{center}
	\hspace{0mm}
    \includegraphics[width=120mm,clip]{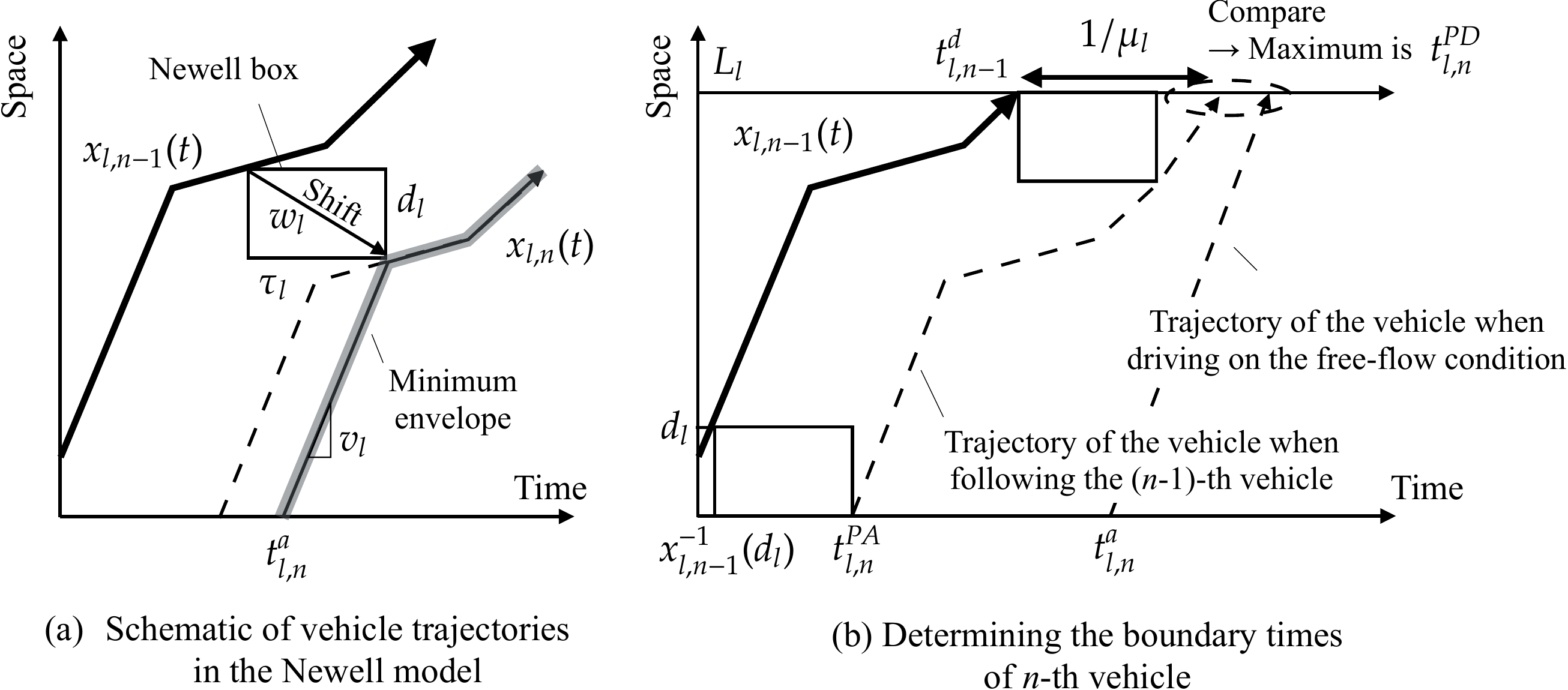}
	\end{center}
    \vspace{-4mm}
	\caption{Schematic of the Newell's car-following model and possible times}
    \vspace{-3mm}
    \label{Fig:NewellBox}
\end{figure}

A link model must satisfy several natural conditions, including the FIFO principle and causality~\citep[][]{Carey2003} to eliminate vehicle movements that are physically implausible. 
We here employ the car-following model by \cite{Newell2002}.
This model describes how the trajectory $x_{l,n}(t)$, which is the position of the $n$-th vehicle (`follower') entering link $l$, depends on the trajectory $x_{l,n-1}(t)$ of the $(n-1)$-th vehicle (`leader'). 
Specifically, the trajectory $x_{l,n}(t)$ is determined such that the spacing and headway are kept at more than the values of the `Newell box' (see Figure~\ref{Fig:NewellBox}(a)).
The Newell box is characterised by two variables: namely, reaction time $\tau_{l}$ and jam spacing $d_{l}$, which can be calculated by
\begin{align*}
\tau_{l} = \cfrac{1}{w_{l}\kappa_{l}},\quad d_{l} = \cfrac{1}{\kappa_{l}},\quad \text{where}\quad \kappa_{l} = \cfrac{(v_{l} + w_{l})q_{l}}{v_{l}w_{l}}.
\end{align*}
Then, for a given trajectory of the leader $x_{l,n-1}(t)$ and link arrival time of the follower $t^{a}_{l,n}$, the trajectory $x_{l,n}(t)$ is calculated as follows: 
\begin{align}
x_{l,n}(t) = \min.\left\{v_{l}(t-t^{a}_{l,n}), x_{l,n-1}(t-\tau_{l}) - d_{l}   \right\}.
\end{align}
The left hand side of the minimum condition means that the follower drives in the free-flow condition; the right hand side means that the follower drives in the congested condition. 

Based on the model, the possible departure/arrival times of the $n$-th vehicle on a link are derived from the trajectory of the $(n-1)$-th vehicle (see Figure~\ref{Fig:NewellBox}(b)).
The earliest possible arrival time of the $n$-th vehicle $t_{l,n}^{PA}$ is determined from $x^{-1}_{l,n-1}(d_{l})$ representing the time that the $(n-1)$-th vehicle arrives at the position $d_{l}$, as follows\footnote{Note that in case that the $(n-1)$-th vehicle stops at position $d_{l}$ and $x^{-1}_{l,n-1}(d_{l})$ has multiple values, the maximum of $x^{-1}_{l,n-1}(d_{l})$ is utilised to determine the possible arrival time.}: 
\begin{align}
t_{l,n}^{PA} = x^{-1}_{l,n-1}(d_{l}) + \tau_{l}.
\end{align}
The earliest departure time $t_{l,n}^{PD}$ is determined from the departure time of the $(n-1)$-th vehicle $t_{l,n-1}^{d}$ and the (actual) arrival time of the $n$-th vehicle $t^{a}_{l,n}$:
\begin{align}
t_{l,n}^{PD} = \max.\left\{ t_{l,n-1}^{d} + \cfrac{1}{\mu_{l}}, t^{a}_{l,n} +\cfrac{L_{l}}{v_{l}}    \right\}. \label{Eq:PossibleDemand}
\end{align}

\subsection{Node model}
A node model determines the actual departure/arrival times of vehicles based on possible times derived from the link model.
Each node belongs to one of the following node types: a normal node that connects one upstream link and one downstream link; a diverge node that connects one upstream link and multiple downstream links; a merge node that connects multiple upstream links and one downstream link; an intersection node that connects multiple upstream and downstream links.
Especially, for nodes other than the normal nodes, the departure/arrival times are determined while taking into account diverging or/and merging behaviour of vehicles.

\begin{figure}[t]
	\begin{center}
	\hspace{0mm}
    \includegraphics[width=70mm,clip]{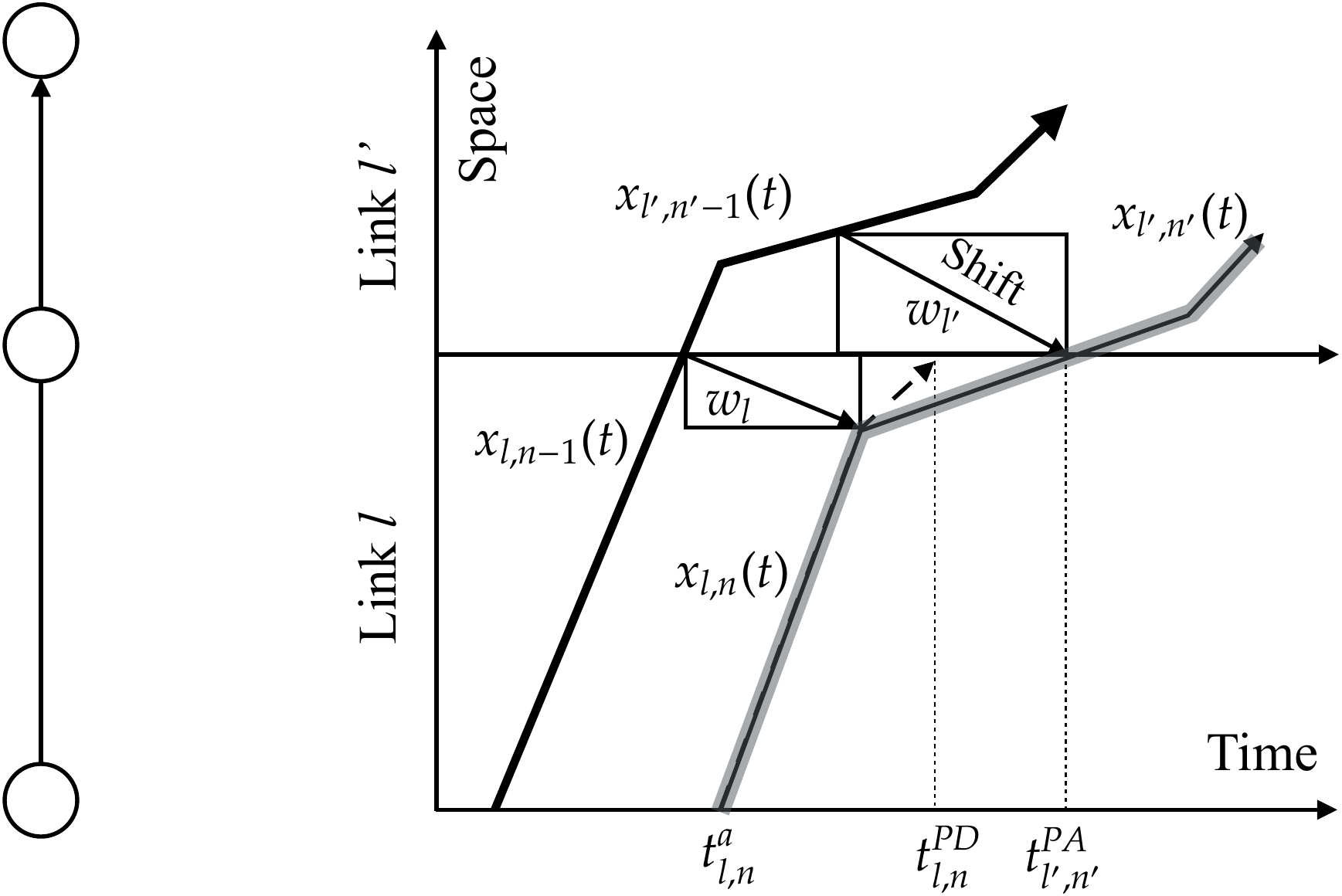}
	\end{center}
    \vspace{-4mm}
	\caption{Trajectories of vehicles passing through a node with one upstream and downstream links}
    \vspace{-3mm}
    \label{Fig:Trajectory_PassNode}
\end{figure}

On a normal node, the actual departure time from the upstream link is determined as the maximum of the possible departure time from the upstream link and the possible arrival time at the downstream link (see Figure~\ref{Fig:Trajectory_PassNode}).
Specifically, when the $n$-th vehicle on the upstream link $l$ enters the downstream link $l'$ as the $n'$-th vehicle, the departure time of the vehicle from link $l$ (i.e. the arrival time at link $l'$) is determined as follows: 
\begin{align}
t^{d}_{l,n} = t^{a}_{l',n'} = \max\left\{ t_{l,n}^{PD}, t_{l',n'}^{PA}   \right\}.\label{Eq:NodeModel}
\end{align}


On a diverge node, we assume that the movement of a follower vehicle is restricted by the leader vehicle on the link which the follower vehicle is going to enter.
Thus, we determine the actual departure/arrival times of a vehicle from Eq.~\eqref{Eq:NodeModel} by regarding the link that is included in the route of the vehicle as $l'$.

On a merge node, we have to determine which vehicle can actually enter the downstream link when more than one vehicle can simultaneously enter the link. 
We here assume that the vehicle on the upstream link having the highest `priority' can depart from the link and enter the downstream link.
The priority at time $t$ of each upstream link of a node is determined according to traffic flow passing through the node before time $t$.
Specifically, the priority is determined such that the vehicles enter from the upstream links in a ratio equal to the ratio of their capacities, e.g. when there exist two upstream links $l_{a}$ and $l_{b}$ of which the ratio of the capacities is $2:1$ and one vehicle has already departed from each link before time $t$, we prioritise the link $l_{a}$\footnote{This can be regarded as an atomic version of the Daganzo's (fluid) merge model~\citep[][]{Daganzo1995}.}.
Then, actual departure/arrival times can be determined by utilising the possible departure time of the prioritised link and possible arrival time of the downstream link from Eq.~\eqref{Eq:NodeModel}.
Note that, on an intersection node, we can determine actual departure/arrival times with a combination of the proposed diverge and merge rules.

\section{Proof of \textbf{Proposition~\ref{Theo:LogitPara}}}\label{App:LogitPara}
This proposition is derived by employing Theorem~6 in \cite{Tatarenko2014}.
Prior to the introduction of the theorem, we show the following notions utilised in the theorem: 

\begin{defi}\textbf{(Scrambling matrix).}
A matrix $\mathbf{P}$ is called as scrambling matrix, if for any two rows, say $\alpha$ and $\beta$, there exists at least one column, $\gamma$, such that both $p_{\alpha\gamma}>0$ and $p_{\beta\gamma}>0$.
\end{defi}




By utilising the notion, the following theorem is proposed: 
\begin{theo}\textbf{(Theorem~6 of \cite{Tatarenko2014}).}
Consider a potential game, and let $\mathbf{B}$ be a pattern matrix of the logit response dynamics in the game; the pattern matrix is the transition probability matrix of the logit response dynamics between two states where the parameter $\beta = 0$.
Let also integer $n<\infty$ be a minimal integer such that the $n$-th power of the matrix $\mathbf{B}$ (i.e. $\mathbf{B}^{n}$) is scrambling.
Then, the logit response dynamics with $\epsilon(\tau) = \ln(\tau+1)/n$ applied to the game guarantees the probabilistic convergence to joint actions to the maximisers of potential function $\Pi$, namely,
\begin{align}
\lim_{\tau\rightarrow \infty} \Pr \{   \mathbf{r}^{\tau} \in \{  \mathbf{r}^{*}  \mid     \Pi(\mathbf{r}^{*}) = \max_{\mathbf{r}} \Pi(\mathbf{r})   \}  \} = 1.
\end{align}

\end{theo}

Thus, it is sufficient for us to prove that the integer $n$ is $ \lceil |\mathcal{P}|/2 \rceil$ in a DSO game with $|\mathcal{P}|$ users.
To prove this, we utilise the following proposition: 
\begin{prop}
Consider a pattern matrix $\mathbf{B}$ of the logit response dynamics in a DSO game.
Consider also two rows $\alpha$ and $\beta$, and corresponding route profiles $\mathbf{r}^{\alpha}$ and $\mathbf{r}^{\beta}$.
Let $\mathcal{P}_{\alpha\beta}$ be the set of users taking different routes in these route profiles, that is, 
\begin{align}
\begin{cases}
r^{\alpha}_{i}\neq r^{\beta}_{i},         &\quad i\in\mathcal{P}_{\alpha\beta}\\
r^{\alpha}_{i}= r^{\beta}_{i} = r_{i},	&\quad i\in\mathcal{P}\setminus \mathcal{P}_{\alpha\beta}
\end{cases}.
\end{align}

\noindent Then, $n = \lceil |\mathcal{P}_{\alpha\beta}|/2 \rceil $ is a minimal integer such that for the two rows $\alpha$ and $\beta$ in $\mathbf{B}^{n}$, there exists at least one column $\gamma$ satisfying the following conditions: 
\begin{align}
b^{ \lceil |\mathcal{P}_{\alpha\beta}|/2 \rceil}_{\alpha \gamma} > 0  \text{ \textcolor{black}{and} } b^{ \lceil |\mathcal{P}_{\alpha\beta}|/2 \rceil}_{\beta \gamma} > 0.
\end{align}
\end{prop}

\begin{figure}[t]
	\begin{center}
	\hspace{0mm}
    \includegraphics[width=80mm,clip]{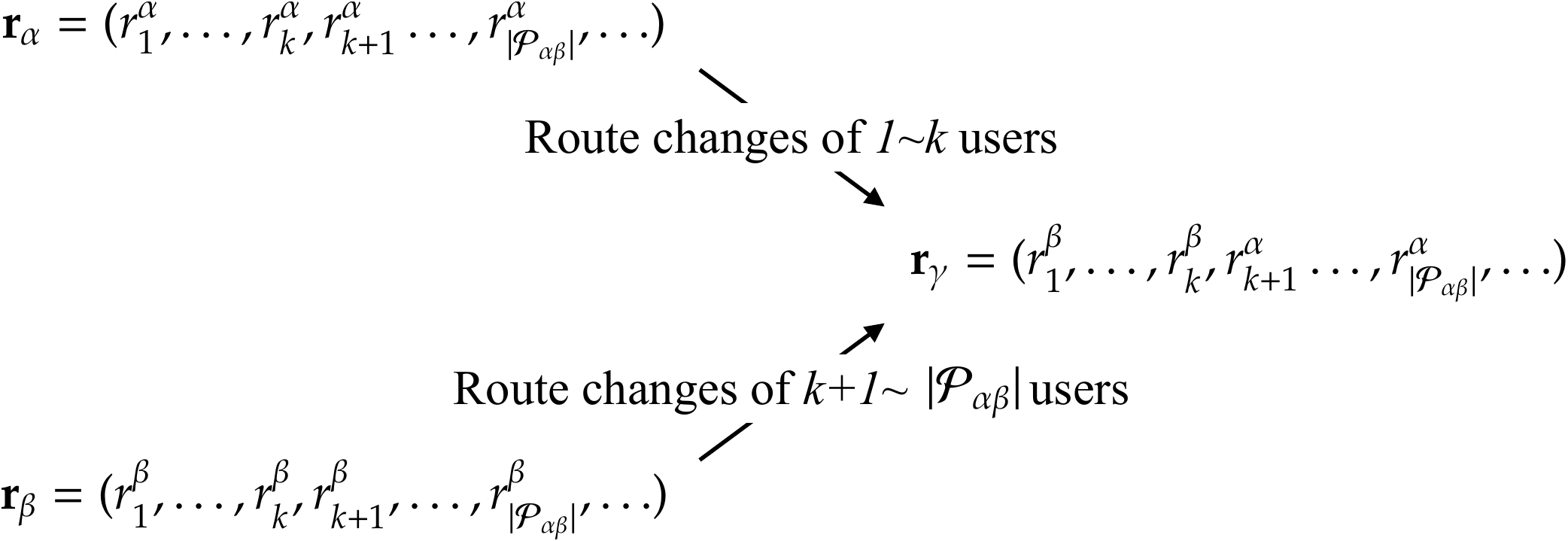}
	\end{center}
    \vspace{-4mm}
	\caption{Schematic of the transition from the route profile $\mathbf{r}^{\alpha}$ and $\mathbf{r}^{\beta}$ to $\mathbf{r}^{\gamma}$}
    \vspace{-3mm}
    \label{Fig:Scrambling}
\end{figure}

\begin{prf}
First of all, each element $(a,b)$ of the matrix $\mathbf{B}^{n}$ can be represented by utilising the corresponding route profiles $\mathbf{r}^{a}$ and $\mathbf{r}^{b}$ as follows: 
\begin{align}
b^{n}_{a b}=
\begin{cases}
>0  &  \text{if $\mathbf{r}^{b}$ is reachable from $\mathbf{r}^{a}$ by route changes of at most $n$ users} \\
0   & \text{otherwise}
\end{cases}.
\end{align}

\noindent This implies that if a route profile $\mathbf{r}^{\gamma}$ is reachable from $\mathbf{r}^{\alpha}$ by route changes of $k$ $(<|\mathcal{P}_{\alpha\beta}|)$ users, at least $|\mathcal{P}_{\alpha\beta}|-k$ users' route changes are required to reach $\mathbf{r}^{\gamma}$ from $\mathbf{r}^{\beta}$.
Specifically, if $\mathbf{r}^{\gamma}$ is a route profile that is reached from $\mathbf{r}^{\alpha}$ by route changes of $k$ users in the set $|\mathcal{P}_{\alpha\beta}|$ to the same routes in $\mathbf{r}^{\beta}$, $\mathbf{r}^{\gamma}$ can be reached from $\mathbf{r}^{\beta}$ by route changes of the other $|\mathcal{P}_{\alpha\beta}|-k$ users to the same routes in $\mathbf{r}^{\alpha}$ (see Figure~\ref{Fig:Scrambling}). 
It follows that by iterating the logit response dynamics $\max\{ k, |\mathcal{P}_{\alpha\beta}|-k  \}$ times, there exists at least one route profile $\mathbf{r}^{\gamma}$ that is reachable from $\mathbf{r}^{\alpha}$ and $\mathbf{r}^{\beta}$, that is, 
\begin{align}
b^{\max\{ k, |\mathcal{P}_{\alpha\beta}|-k  \}  }_{\mathbf{r}^{\alpha} \mathbf{r}^{\gamma}} > 0  \text{ \textcolor{black}{and} }b^{\max\{ k, |\mathcal{P}_{\alpha\beta}|-k  \}}_{\mathbf{r}^{\beta} \mathbf{r}^{\gamma}} > 0,
\end{align}
For the given $\mathcal{P}_{\alpha\beta}$, it is obvious that the max function is minimised by $k= \lceil |\mathcal{P}_{\alpha\beta}|/2 \rceil$.
Thus, the proposition is proved.\qed

\end{prf}

\noindent From this proposition, in a DSO game with $|\mathcal{P}|$, it follows that if $n$ is equal to or more than $\lceil |\mathcal{P}|/2 \rceil$, $\mathbf{B}^{n}$ becomes a scrambling matrix since the maximum number of users who take different routes is $|\mathcal{P}|$.
From this fact and the theorem proposed by \cite{Tatarenko2014}, the proposition is derived.\qed

\section{Proof of \textbf{Corollary~\ref{Coro:FixedTollAdjust}}}\label{App:FixedTollAdjust}

We prove the corollary by contradiction.
Suppose that there exists another route profile $\mathbf{r}^{**}$ which is Nash equilibrium.
Among the users taking different routes in these route profiles, we consider the user $i$ departing earliest.  
From \textbf{Proposition~\ref{Prop:OrderingProperty}} and the fact that the users departing earlier than user $i$ take the same routes, the following relationship holds: 
\begin{align}
U^{F}_{i}(r, \mathbf{r}^{*}_{-i}) = U^{F}_{i}(r, \mathbf{r}^{**}_{-i}),\quad \forall r\in\mathcal{R}_{i}.
\end{align}
However, by combining this relationship with Eq.~\eqref{Eq:TollAdjustment}, we obtain the following relationship:
\begin{align}
U^{F}_{i}(r^{**}_{i}, \mathbf{r}^{**}_{-i}) = U^{F}_{i}(r^{**}_{i}, \mathbf{r}^{*}_{-i})<   U^{F}_{i}(r^{*}_{i}, \mathbf{r}^{*}_{-i})= U^{F}_{i}(r^{*}_{i}, \mathbf{r}^{**}_{-i}) .
\end{align}
This contradicts the statement that $\mathbf{r}^{**}$ is a Nash equilibrium, i.e. there do not exist multiple equilibrium states.
Thus, the proposition is proved.

\section{Proof of \textbf{Theorem~\ref{Theo:DUE-FCPisWAG}}}\label{App:DUE-FCPisWAG}

We prove the existence of a better response path ending at a Nash equilibrium state from every route profile in a constructive manner: we construct an algorithm that is guaranteed to move any initial traffic state with route profile $\mathbf{r}$ into an equilibrium state $\mathbf{r}^{*}$ by a sequence of better responses.
In the algorithm, the users are divided into two groups: (i) Group A contains the users who are currently on their \textit{ex-post} best response routes , and (ii) Group B contains users who are not.
The sets of users in Group A and Group B are denoted by $\mathcal{P}_{A}$ and $\mathcal{P}_{B}$, respectively.
We also denote by $\mathbf{r}_{A}$ and $\mathbf{r}_{B}$ the route profiles of users in Group A and Group B, respectively (i.e. $\mathbf{r} = (\mathbf{r}_{A}, \mathbf{r}_{B})$).

We then propose the following algorithm: 

\vspace{2mm}
\begin{algorithm}[H]\label{Algo:PrfWAG}
\caption{\textbf{Leading a traffic state to a Nash equilibrium state by better responses}}

{
0. \textbf{Initialisation}: Set $m = 0$ where $m$ represents the iteration counter. 
Initialise all users to be in Group B, i.e. $(\mathcal{P}_{A}^{m}, \mathcal{P}_{B}^{m}) = (\emptyset, \mathcal{P})$ and $\mathbf{r}^{m} = (\mathbf{r}^{m}_{A}, \mathbf{r}^{m}_{B}) = (\mathbf{0}, \mathbf{r})$.

1. \textbf{Pick up the user departing from the origin earliest}: 
Pick up the vehicle departing from the origin earliest from Group B.
The vehicle is denoted by $i^{m} \in \mathcal{P}_{B}^{m}$.
Note that such a vehicle can uniquely be chosen for each iteration because we assume that all users with the same origin depart at different times.

2. \textbf{Update the route profile}: 
Find a best response route of user $i^{m}$ to the current route profile $\mathbf{r}^{m}$, i.e. find a route $r^{*}\in\mathcal{R}_{i^{m}}$ satisfying the following condition: $U_{i^{m}}^{F}(r, \mathbf{r}^{m}_{-i^{m}}) \leq U_{i^{m}}^{F}(r^{*}, \mathbf{r}^{m}_{-i^{m}}), \ \forall r \in \mathcal{R}_{i^{m}}.$
Compare the utility when user $i^{m}$ takes the current route $r^{m}_{i^{m}}$ and the utility when taking route $r^{*}$.
Then, 

\begin{description}
\item[(a)] If the former utility is lower than the latter (i.e. $U^{F}_{i^{m}}(r^{m}_{i^{m}}, \mathbf{r}^{m}_{-i^{m}}) < U^{F}_{i^{m}}(r^{*}, \mathbf{r}^{m}_{-i^{m}})$), change the route of user $i^{m}$ to $r^{*}$, i.e. $r_{i^{m}}^{m+1}:= r^{*}$.
Note that the routes of the other users remain the same: $r_{j}^{m+1}:= r_{j}^{m}$ for all $j\in\mathcal{P}\setminus \{ i^{m} \}$.
\item[(b)] If both are the same (i.e. $U^{F}_{i^{m}}(r^{m}_{i^{m}}, \mathbf{r}^{m}_{-i^{m}}) = U^{F}_{i^{m}}(r^{*}, \mathbf{r}^{m}_{-i^{m}})$), the routes of all users are not changed in the current iteration, i.e. $r_{j}^{m+1}:= r_{j}^{m}$ for all $j\in\mathcal{P}$.
\end{description}

3. \textbf{Update the sets of users and judge the convergence}: Let $\mathcal{P}_{A}^{m+1}:=\mathcal{P}_{A}^{m} \cup \{i^{m}\}$ and  $\mathcal{P}_{B}^{m+1} := \mathcal{P}_{B}^{m} \setminus \{i^{m}\}$.
If $\mathcal{P}_{B}^{m+1} \neq \emptyset$, let $m:= m+1$ and go back Step 1.
If $\mathcal{P}_{B}^{m+1} = \emptyset$, then terminate the algorithm; ${\bf r}^{m+1} = {\bf r}^{m+1}_{A}$ is a Nash equilibrium state. 
}
\end{algorithm}
\vspace{2mm}

\noindent This algorithm derives an equilibrium state by assigning users to their best response routes in the order of their departure times (please refer to \cite{Satsukawa2019} for the details of the role of each step).
Thanks to \textbf{Proposition~\ref{Prop:OrderingProperty}}, it is guaranteed that the utility of the user in Group A is not changed by route change of any user in Group B. 
Therefore, once a user is transferred to Group A by Step 2, i.e. the user takes a best response route, the route remains a best response route for the user regardless of route choices of the users in Group B; the route is an ex-post best response route for the user.
It follows that this algorithm can find a Nash equilibrium state in exactly $|\mathcal{P}|$ iterations.
Therefore, there exists a better response path ending at a Nash equilibrium state from an arbitrary initial route profile in a DUE-FCP game in an SBPR-1 network.	\qed

\bibliography{ISTTT24}

\begin{thebibliography}{40}
\expandafter\ifx\csname natexlab\endcsname\relax\def\natexlab#1{#1}\fi
\providecommand{\url}[1]{\texttt{#1}}
\providecommand{\href}[2]{#2}
\providecommand{\path}[1]{#1}
\providecommand{\DOIprefix}{doi:}
\providecommand{\ArXivprefix}{arXiv:}
\providecommand{\URLprefix}{URL: }
\providecommand{\Pubmedprefix}{pmid:}
\providecommand{\doi}[1]{\href{http://dx.doi.org/#1}{\path{#1}}}
\providecommand{\Pubmed}[1]{\href{pmid:#1}{\path{#1}}}
\providecommand{\bibinfo}[2]{#2}
\ifx\xfnm\relax \def\xfnm[#1]{\unskip,\space#1}\fi
\bibitem[{Akamatsu et~al.(2015)Akamatsu, Wada and Hayashi}]{Akamatsu2015}
\bibinfo{author}{Akamatsu, T.}, \bibinfo{author}{Wada, K.},
  \bibinfo{author}{Hayashi, S.}, \bibinfo{year}{2015}.
\newblock \bibinfo{title}{{The corridor problem with discrete multiple
  bottlenecks}}.
\newblock \bibinfo{journal}{Transportation Research Part B: Methodological}
  \bibinfo{volume}{81}, \bibinfo{pages}{808--829}.
\bibitem[{Blume(1993)}]{Blume1993}
\bibinfo{author}{Blume, L.E.}, \bibinfo{year}{1993}.
\newblock \bibinfo{title}{{The statistical mechanics of strategic
  interaction}}.
\newblock \bibinfo{journal}{Games and Economic Behavior} \bibinfo{volume}{5},
  \bibinfo{pages}{387--424}.
\bibitem[{Carey(1992)}]{Carey1992}
\bibinfo{author}{Carey, M.}, \bibinfo{year}{1992}.
\newblock \bibinfo{title}{{Nonconvexity of the dynamic traffic assignment
  problem}}.
\newblock \bibinfo{journal}{Transportation Research Part B}
  \bibinfo{volume}{26}, \bibinfo{pages}{127--133}.
\bibitem[{Carey et~al.(2003)Carey, Ge and McCartney}]{Carey2003}
\bibinfo{author}{Carey, M.}, \bibinfo{author}{Ge, Y.E.},
  \bibinfo{author}{McCartney, M.}, \bibinfo{year}{2003}.
\newblock \bibinfo{title}{{A whole-link travel-time model with desirable
  properties}}.
\newblock \bibinfo{journal}{Transportation Science} \bibinfo{volume}{37},
  \bibinfo{pages}{83--96}.
\bibitem[{Carey and Srinivasan(1993)}]{Carey1993}
\bibinfo{author}{Carey, M.}, \bibinfo{author}{Srinivasan, A.},
  \bibinfo{year}{1993}.
\newblock \bibinfo{title}{{Externalities, average and marginal costs, and tolls
  on congested networks with time-varying flows}}.
\newblock \bibinfo{journal}{Operations Research} \bibinfo{volume}{41},
  \bibinfo{pages}{217--231}.
\bibitem[{Carey and Watling(2012)}]{Carey2012}
\bibinfo{author}{Carey, M.}, \bibinfo{author}{Watling, D.},
  \bibinfo{year}{2012}.
\newblock \bibinfo{title}{{Dynamic traffic assignment approximating the
  kinematic wave model: System optimum, marginal costs, externalities and
  tolls}}.
\newblock \bibinfo{journal}{Transportation Research Part B: Methodological}
  \bibinfo{volume}{46}, \bibinfo{pages}{634--648}.
\bibitem[{Daganzo(1994)}]{Daganzo1994}
\bibinfo{author}{Daganzo, C.F.}, \bibinfo{year}{1994}.
\newblock \bibinfo{title}{{The cell transmission model: A dynamic
  representation of highway traffic consistent with the hydrodynamic theory}}.
\newblock \bibinfo{journal}{Transportation Research Part B: Methodological}
  \bibinfo{volume}{28}, \bibinfo{pages}{269--287}.
\bibitem[{Daganzo(1995)}]{Daganzo1995}
\bibinfo{author}{Daganzo, C.F.}, \bibinfo{year}{1995}.
\newblock \bibinfo{title}{{The cell transmission model, part II: Network
  traffic}}.
\newblock \bibinfo{journal}{Transportation Research Part B: Methodological}
  \bibinfo{volume}{29}, \bibinfo{pages}{79--93}.
\bibitem[{Garcia et~al.(2000)Garcia, Reaume and Smith}]{Garcia2000}
\bibinfo{author}{Garcia, A.}, \bibinfo{author}{Reaume, D.},
  \bibinfo{author}{Smith, R.L.}, \bibinfo{year}{2000}.
\newblock \bibinfo{title}{{Fictitious play for finding system optimal routings
  in dynamic traffic networks}}.
\newblock \bibinfo{journal}{Transportation Research Part B: Methodological}
  \bibinfo{volume}{34}, \bibinfo{pages}{147--156}.
\bibitem[{Ghali and Smith(1995)}]{Ghali1995}
\bibinfo{author}{Ghali, M.}, \bibinfo{author}{Smith, M.}, \bibinfo{year}{1995}.
\newblock \bibinfo{title}{{A model for the dynamic system optimum traffic
  assignment problem}}.
\newblock \bibinfo{journal}{Transportation Research Part B: Methodological}
  \bibinfo{volume}{29}, \bibinfo{pages}{155--170}.
\bibitem[{Kuwahara and Akamatsu(1993)}]{Kuwahara1993}
\bibinfo{author}{Kuwahara, M.}, \bibinfo{author}{Akamatsu, T.},
  \bibinfo{year}{1993}.
\newblock \bibinfo{title}{{Dynamic equilibrium assignment with queues for a
  one-to-many OD pattern}}, in: \bibinfo{editor}{Daganzo, C.F.} (Ed.),
  \bibinfo{booktitle}{Proceedings of the 12th International Symposium on the
  Theory of Traffic Flow and Transportation}, \bibinfo{publisher}{Elsevier},
  \bibinfo{address}{Berkeley}. pp. \bibinfo{pages}{185--204}.
\bibitem[{Kuwahara et~al.(2001)Kuwahara, Yoshii and Kumagai}]{Kuwahara2001}
\bibinfo{author}{Kuwahara, M.}, \bibinfo{author}{Yoshii, T.},
  \bibinfo{author}{Kumagai, K.}, \bibinfo{year}{2001}.
\newblock \bibinfo{title}{{An analysis on dynamic system optimal assignment and
  ramp control on a simple network}}.
\newblock \bibinfo{journal}{Journal of Japan Society of Civil Engineers}
  \bibinfo{volume}{2001}, \bibinfo{pages}{59--71 [In Japanese]}.
\bibitem[{Marden and Shamma(2012)}]{Marden2012b}
\bibinfo{author}{Marden, J.R.}, \bibinfo{author}{Shamma, J.S.},
  \bibinfo{year}{2012}.
\newblock \bibinfo{title}{{Revisiting log-linear learning: Asynchrony,
  completeness and payoff-based implementation}}.
\newblock \bibinfo{journal}{Games and Economic Behavior} \bibinfo{volume}{75},
  \bibinfo{pages}{788--808}.
\bibitem[{Merchant and Nemhauser(1978a)}]{Merchant1978}
\bibinfo{author}{Merchant, D.K.}, \bibinfo{author}{Nemhauser, G.L.},
  \bibinfo{year}{1978}a.
\newblock \bibinfo{title}{{A model and an algorithm for the dynamic traffic
  assignment problems}}.
\newblock \bibinfo{journal}{Transportation Science} \bibinfo{volume}{12},
  \bibinfo{pages}{183--199}.
\bibitem[{Merchant and Nemhauser(1978b)}]{Merchant1978a}
\bibinfo{author}{Merchant, D.K.}, \bibinfo{author}{Nemhauser, G.L.},
  \bibinfo{year}{1978}b.
\newblock \bibinfo{title}{{Optimality conditions for a dynamic traffic
  assignment model}}.
\newblock \bibinfo{journal}{Transportation Science} \bibinfo{volume}{12},
  \bibinfo{pages}{200--207}.
\bibitem[{Meyn and Tweedie(2009)}]{Meyn2009}
\bibinfo{author}{Meyn, S.}, \bibinfo{author}{Tweedie, R.L.},
  \bibinfo{year}{2009}.
\newblock \bibinfo{title}{Markov Chains and Stochastic Stability}.
\newblock \bibinfo{edition}{2nd} ed., \bibinfo{publisher}{Cambridge University
  Press}, \bibinfo{address}{USA}.
\bibitem[{Monderer and Shapley(1996)}]{Monderer1996}
\bibinfo{author}{Monderer, D.}, \bibinfo{author}{Shapley, L.S.},
  \bibinfo{year}{1996}.
\newblock \bibinfo{title}{{Potential games}}.
\newblock \bibinfo{journal}{Games and Economic Behavior} \bibinfo{volume}{14},
  \bibinfo{pages}{124--143}.
\bibitem[{Mu{\~{n}}oz and Laval(2006)}]{Munoz2006}
\bibinfo{author}{Mu{\~{n}}oz, J.C.}, \bibinfo{author}{Laval, J.A.},
  \bibinfo{year}{2006}.
\newblock \bibinfo{title}{{System optimum dynamic traffic assignment graphical
  solution method for a congested freeway and one destination}}.
\newblock \bibinfo{journal}{Transportation Research Part B: Methodological}
  \bibinfo{volume}{40}, \bibinfo{pages}{1--15}.
\bibitem[{Newell(2002)}]{Newell2002}
\bibinfo{author}{Newell, G.F.}, \bibinfo{year}{2002}.
\newblock \bibinfo{title}{{A simplified car-following theory: A lower order
  model}}.
\newblock \bibinfo{journal}{Transportation Research Part B: Methodological}
  \bibinfo{volume}{36}, \bibinfo{pages}{195--205}.
\bibitem[{Nie(2011)}]{Nie2011}
\bibinfo{author}{Nie, Y.M.}, \bibinfo{year}{2011}.
\newblock \bibinfo{title}{{A cell-based Merchant-Nemhauser model for the system
  optimum dynamic traffic assignment problem}}.
\newblock \bibinfo{journal}{Transportation Research Part B: Methodological}
  \bibinfo{volume}{45}, \bibinfo{pages}{329--342}.
\bibitem[{Qian et~al.(2012)Qian, Shen, Zhang, Sean, Shen and Zhang}]{Qian2012}
\bibinfo{author}{Qian, Z.S.}, \bibinfo{author}{Shen, W.},
  \bibinfo{author}{Zhang, H.M.}, \bibinfo{author}{Sean, Z.},
  \bibinfo{author}{Shen, W.}, \bibinfo{author}{Zhang, H.M.},
  \bibinfo{year}{2012}.
\newblock \bibinfo{title}{{System-optimal dynamic traffic assignment with and
  without queue spillback: Its path-based formulation and solution via
  approximate path marginal cost}}.
\newblock \bibinfo{journal}{Transportation Research Part B: Methodological}
  \bibinfo{volume}{46}, \bibinfo{pages}{874--893}.
\bibitem[{Sandholm(2002)}]{Sandholm2002}
\bibinfo{author}{Sandholm, W.H.}, \bibinfo{year}{2002}.
\newblock \bibinfo{title}{{Evolutionary implementation and congestion
  pricing}}.
\newblock \bibinfo{journal}{Review of Economic Studies} \bibinfo{volume}{69},
  \bibinfo{pages}{667--689}.
\bibitem[{Sandholm(2005)}]{Sandholm2005}
\bibinfo{author}{Sandholm, W.H.}, \bibinfo{year}{2005}.
\newblock \bibinfo{title}{{Negative externalities and evolutionary
  implementation}}.
\newblock \bibinfo{journal}{Review of Economic Studies} \bibinfo{volume}{72},
  \bibinfo{pages}{885--915}.
\bibitem[{Sandholm(2007)}]{Sandholm2007}
\bibinfo{author}{Sandholm, W.H.}, \bibinfo{year}{2007}.
\newblock \bibinfo{title}{{Pigouvian pricing and stochastic evolutionary
  implementation}}.
\newblock \bibinfo{journal}{Journal of Economic Theory} \bibinfo{volume}{132},
  \bibinfo{pages}{367--382}.
\bibitem[{Satsukawa et~al.(2019)Satsukawa, Wada and Iryo}]{Satsukawa2019}
\bibinfo{author}{Satsukawa, K.}, \bibinfo{author}{Wada, K.},
  \bibinfo{author}{Iryo, T.}, \bibinfo{year}{2019}.
\newblock \bibinfo{title}{{Stochastic stability of dynamic user equilibrium in
  unidirectional networks: Weakly acyclic game approach}}.
\newblock \bibinfo{journal}{Transportation Research Part B: Methodological}
  \bibinfo{volume}{125}, \bibinfo{pages}{229--247}.
\bibitem[{Shen et~al.(2007)Shen, Nie and Zhang}]{Shen2007}
\bibinfo{author}{Shen, W.}, \bibinfo{author}{Nie, Y.}, \bibinfo{author}{Zhang,
  H.M.}, \bibinfo{year}{2007}.
\newblock \bibinfo{title}{{On path marginal cost analysis and its relation to
  dynamic system-optimal traffic assignment}}, in: \bibinfo{editor}{Allsop,
  R.E.}, \bibinfo{editor}{Bell, M.G.H.}, \bibinfo{editor}{Heydecker, B.}
  (Eds.), \bibinfo{booktitle}{Proceedings of the 17th international symposium
  on transportation and traffic theory}, \bibinfo{publisher}{Elsevier},
  \bibinfo{address}{London , England}. pp. \bibinfo{pages}{327--360}.
\bibitem[{Shen and Zhang(2009)}]{Shen2009}
\bibinfo{author}{Shen, W.}, \bibinfo{author}{Zhang, H.M.},
  \bibinfo{year}{2009}.
\newblock \bibinfo{title}{{On the morning commute problem in a corridor network
  with multiple bottlenecks: Its system-optimal traffic flow patterns and the
  realizing tolling scheme}}.
\newblock \bibinfo{journal}{Transportation Research Part B: Methodological}
  \bibinfo{volume}{43}, \bibinfo{pages}{267--284}.
\bibitem[{Shen and Zhang(2014)}]{Shen2014}
\bibinfo{author}{Shen, W.}, \bibinfo{author}{Zhang, H.M.},
  \bibinfo{year}{2014}.
\newblock \bibinfo{title}{{System optimal dynamic traffic assignment:
  Properties and solution procedures in the case of a many-to-one network}}.
\newblock \bibinfo{journal}{Transportation Research Part B: Methodological}
  \bibinfo{volume}{65}, \bibinfo{pages}{1--17}.
\bibitem[{Smits et~al.(2015)Smits, Bliemer, Pel and van Arem}]{Smits2015}
\bibinfo{author}{Smits, E.S.}, \bibinfo{author}{Bliemer, M.C.},
  \bibinfo{author}{Pel, A.J.}, \bibinfo{author}{van Arem, B.},
  \bibinfo{year}{2015}.
\newblock \bibinfo{title}{{A family of macroscopic node models}}.
\newblock \bibinfo{journal}{Transportation Research Part B: Methodological}
  \bibinfo{volume}{74}, \bibinfo{pages}{20--39}.
\bibitem[{Stewart(2009)}]{Stewart2009}
\bibinfo{author}{Stewart, W.J.}, \bibinfo{year}{2009}.
\newblock \bibinfo{title}{Probability, Markov chains, queues, and simulation}.
\newblock \bibinfo{publisher}{Princeton university press}.
\bibitem[{Tamp{\`{e}}re et~al.(2011)Tamp{\`{e}}re, Corthout, Cattrysse and
  Immers}]{Tampere2011}
\bibinfo{author}{Tamp{\`{e}}re, C.M.}, \bibinfo{author}{Corthout, R.},
  \bibinfo{author}{Cattrysse, D.}, \bibinfo{author}{Immers, L.H.},
  \bibinfo{year}{2011}.
\newblock \bibinfo{title}{{A generic class of first order node models for
  dynamic macroscopic simulation of traffic flows}}.
\newblock \bibinfo{journal}{Transportation Research Part B: Methodological}
  \bibinfo{volume}{45}, \bibinfo{pages}{289--309}.
\bibitem[{Tatarenko(2014)}]{Tatarenko2014}
\bibinfo{author}{Tatarenko, T.}, \bibinfo{year}{2014}.
\newblock \bibinfo{title}{{Proving convergence of log-linear learning in
  potential games}}.
\newblock \bibinfo{journal}{Proceedings of the American Control Conference} ,
  \bibinfo{pages}{972--977}.
\bibitem[{Wada et~al.(2019)Wada, Satsukawa, Smith and Akamatsu}]{Wada2019}
\bibinfo{author}{Wada, K.}, \bibinfo{author}{Satsukawa, K.},
  \bibinfo{author}{Smith, M.J.}, \bibinfo{author}{Akamatsu, T.},
  \bibinfo{year}{2019}.
\newblock \bibinfo{title}{{Network throughput under dynamic user equilibrium :
  Queue spillback, paradox and traffic control}}.
\newblock \bibinfo{journal}{Transportation Research Part B: Methodological}
  \bibinfo{volume}{126}, \bibinfo{pages}{391--413}.
\bibitem[{Young(1993)}]{Young1993}
\bibinfo{author}{Young, H.P.}, \bibinfo{year}{1993}.
\newblock \bibinfo{title}{{The evolution of conventions}}.
\newblock \bibinfo{journal}{Econometrica} \bibinfo{volume}{61},
  \bibinfo{pages}{57--84}.
\bibitem[{Young(2004)}]{Young2004}
\bibinfo{author}{Young, H.P.}, \bibinfo{year}{2004}.
\newblock \bibinfo{title}{{Strategic Learning and Its Limits}}.
\newblock \bibinfo{publisher}{Oxford University Press, USA}.
\bibitem[{Zhang and Shen(2010)}]{Zhang2010}
\bibinfo{author}{Zhang, H.M.}, \bibinfo{author}{Shen, W.},
  \bibinfo{year}{2010}.
\newblock \bibinfo{title}{{Access control policies without inside queues: Their
  properties and public policy implications}}.
\newblock \bibinfo{journal}{Transportation Research Part B: Methodological}
  \bibinfo{volume}{44}, \bibinfo{pages}{1132--1147}.
\bibitem[{Zhang and Qian(2020)}]{Zhang2020}
\bibinfo{author}{Zhang, P.}, \bibinfo{author}{Qian, S.}, \bibinfo{year}{2020}.
\newblock \bibinfo{title}{{Path-based system optimal dynamic traffic assignment
  : A subgradient approach}}.
\newblock \bibinfo{journal}{Transportation Research Part B: Methodological}
  \bibinfo{volume}{134}, \bibinfo{pages}{41--63}.
\bibitem[{Zhao and Leclercq(2018)}]{Zhao2018}
\bibinfo{author}{Zhao, C.L.L.}, \bibinfo{author}{Leclercq, L.},
  \bibinfo{year}{2018}.
\newblock \bibinfo{title}{{Graphical solution for system optimum dynamic
  traffic assignment with day-based incentive routing strategies}}.
\newblock \bibinfo{journal}{Transportation Research Part B: Methodological}
  \bibinfo{volume}{117}, \bibinfo{pages}{87--100}.
\bibitem[{Zhu and Ukkusuri(2013)}]{Zhu2013}
\bibinfo{author}{Zhu, F.}, \bibinfo{author}{Ukkusuri, S.V.},
  \bibinfo{year}{2013}.
\newblock \bibinfo{title}{{A cell based dynamic system optimum model with
  non-holding back flows}}.
\newblock \bibinfo{journal}{Transportation Research Part C: Emerging
  Technologies} \bibinfo{volume}{36}, \bibinfo{pages}{367--380}.
\bibitem[{Ziliaskopoulos(2000)}]{Ziliaskopoulos2000}
\bibinfo{author}{Ziliaskopoulos, A.K.}, \bibinfo{year}{2000}.
\newblock \bibinfo{title}{{A linear programming model for the single
  destination system optimum dynamic traffic assignment problem}}.
\newblock \bibinfo{journal}{Transportation Science} \bibinfo{volume}{34},
  \bibinfo{pages}{36--49}.

\end{thebibliography}
\bibliographystyle{elsarticle-harv}




\end{document}